\documentclass{article}

\usepackage{arxiv}

\usepackage[utf8]{inputenc} % allow utf-8 input
\usepackage[T1]{fontenc}    % use 8-bit T1 fonts
\usepackage{hyperref}       % hyperlinks
\usepackage{url}            % simple URL typesetting
\usepackage{booktabs}       % professional-quality tables
\usepackage{amsfonts}       % blackboard math symbols
\usepackage{amsmath}
\usepackage{nicefrac}       % compact symbols for 1/2, etc.
\usepackage{microtype}      % microtypography
\usepackage{lipsum}
\usepackage{tikz}
\usepackage{subfig}
\usepackage[ruled]{algorithm2e}

\newcommand{\R}{\mathbb{R}}
\renewcommand{\P}{\mathbb{P}}
\newcommand{\F}{\mathcal{F}}
\newcommand{\N}{\mathbb{N}}
\newcommand{\D}{\mathcal{D}}
\newcommand{\E}{\mathbb{E}}
\newcommand{\RR}{\mathcal{R}}
\newcommand{\NN}{\mathcal{N}}
\newcommand{\M}{\mathcal{M}}
\renewcommand{\S}{\mathcal{S}}
\newcommand{\T}{\mathcal{T}}

\newcommand{\veps}{\varepsilon}
\newcommand{\hsigma}{\hat{\sigma}}
\newcommand{\Var}{\text{Var}}
\newcommand{\norm}[1]{\left\Vert #1\right\Vert}
\usepackage{amsthm}
\usepackage{amsfonts}
\newtheorem{theorem}{Theorem}
\newtheorem{proposition}{Proposition}
\usepackage{bbm}

\title{Transfer Learning for Linear Regression: a Statistical Test of Gain}

\author{
  David~Obst\\
  EDF R\&D\\
  Aix-Marseille Université\\
  \texttt{david.obst@edf.fr} \\
   \And
  Badih~Ghattas \\
  Aix-Marseille Université\\
  \texttt{badih.ghattas@univ-amu.fr} \\
  \And
  Jairo~Cugliari \\
  Université Lumière Lyon 2\\
  \texttt{Jairo.Cugliari@univ-lyon2.fr} \\
  \And
  Georges~Oppenheim \\
  Université Paris-Est Marne-la-Vallée\\
  \texttt{georges.oppenheim@gmail.com} \\
  \And
  Sandra~Claudel\\
  EDF R\&D\\
  \texttt{sandra.claudel@edf.fr} \\
  \And
  Yannig~Goude\\
  EDF R\&D\\
  \texttt{yannig.goude@edf.fr} \\
}

\begin{document}
\maketitle

\begin{abstract}
 Transfer learning, also referred as knowledge transfer, aims at reusing knowledge from a source dataset to a similar target one. While many empirical studies illustrate the benefits of transfer learning, few theoretical results are established especially for regression problems. In this paper a theoretical framework for the problem of parameter transfer for the linear model is proposed. It is shown that the quality of transfer for a new input vector $x$ depends on its representation in an eigenbasis involving the parameters of the problem. Furthermore a statistical test is constructed to predict whether a fine-tuned model has a lower prediction quadratic risk than the base target model for an unobserved sample. Efficiency of the test is illustrated on synthetic data as well as real electricity consumption data.
\end{abstract}

% keywords can be removed
\keywords{Linear regression \and Transfer learning \and Statistical test \and Fine-tuning \and Transfer theory}

\section{Introduction}
\label{intro}

%Definir le cadre: peu de donnees T et bcp de S avec un lien
%Definir les types de transfer: param, feature, instance & relation
%Pas besoin de def Inductive, Transductive et Unsupervised.
%Pb: peu de papiers theoriques, surtout en regression.

We consider the situation where we want to perform predictions for a target task $\T$ but for which training data is limited, either because it is difficult to label or only newly available. We have at our disposal another source task $\S$ related (in some sense) to $\T$ for which we have plenty of data. Our goal, which is the one of transfer learning, is to leverage information from $\S$ to improve the task results on $\T$ \cite{weiss2016survey}. For instance we could be interested in performing forecasts for a group of newly arrived customers $\T$ while we have a group of long time customers $\S$. We expect that the relatedness between $\S$ and $\T$ and the information at our disposal for $\S$ helps to improve our forecasting model for the newly arrived customers. Transfer can be performed on multiple levels. \cite{pan2009survey} define four categories of transfer: parameter, instances, features and relationship transfer. In this paper we focus on the first, more specifically in the framework of linear regression. Transfer for the linear model has been studied in the literature on multiple occasions, but the essence of results vary. In \cite{lounici2009taking,lounici2011oracle} the authors consider the setting of sparse multi-task learning in high dimension with a common sparsity pattern within the regression vectors. They obtain oracle inequalities on the prediction error, albeit for the same data on which the parameters were learned. \cite{maurer2006bounds} establishes bounds on the average prediction error over $m$ tasks in the linear classification setting, but do not investigate in which situation learning on multiple sources could be beneficial for a specific target. Furthermore in both aforementioned papers the results remain mainly theoretical, and the setting of multi-task learning is slightly different than the one of transfer. Bayesian approaches, where a prior is constructed with the help of the source data and then the posterior is obtained with the target one are also popular for parameter transfer. For instance in \cite{launay2015construction} it is applied for the problem of electricity load forecasting, with good results. \cite{bouveyron2010adaptive} propose another method to transfer the estimated coefficients of a linear model when the number of target samples is highly limited. After estimating the regression coefficients on the source set, the target vector is obtained by a linear transformation of the previous one with constraints on the transformation. They demonstrated the efficiency of their approach on house pricing data among others, showing significant improvement over learning from scratch on the target set. \cite{chen2015data} suggest another approach of transfer, where the estimator of the target vector is built either by directly minimizing a combination of the quadratic losses of both sets (referred as data pooling in their paper) or by constructing a convex combination of the individual estimators. They proceed to investigate the properties of the estimators, with one of their main results being the optimality of the convex combination under certain conditions. \cite{dar2020double} recently studied transfer of a set of parameters for linear regression in the restricted setting of gaussian and independent features. They showed the existence of a phenomenon of "double double descent", with transfer being beneficial under certain conditions of under or over-parametrization of the tasks. However in the aforementioned papers one common practice in transfer learning is missing, namely fine-tuning. It has shown a lot of success in recent years, notably for neural networks and allows for more flexibility than just combining estimators. It consists in reusing a part of the learned parameters on the source (for instance neural network layers) and adjusting them on the target with a few gradient iterations \cite{shin2016deep}. Furthermore the problem of negative transfer, i.e. when the transfer procedure may be detrimental, is not fully adressed in the aforementioned papers especially for new observations. \cite{fawaz2018transfer} approach the problem empirically: after defining a distance between datasets based on the dynamic time warping (DTW), they show that in general negative transfer will happen when the defined distance between source and target is large. \cite{ben2010theory} indirectly address the issue of negative transfer for the problem of binary classification. Considering the transfer problem as a special case of a multi-task objective, not only do they obtain an upper bound on the transfer prediction error, but they also prove the existence of phases depending on the number of samples $N_S$ and $N_T$ available for source and target respectively. Finally a domination inequality is established in \cite{chen2015data} for their optimal estimator, i.e. meaning that in certain cases transfer is bound to be beneficial. However it is valid only under specific assumptions on the observations that are usually not met in practice.
Therefore we define the problem of negative transfer  as following. We have at our disposal two parametric models, $\M_\T$ of estimated parameter $\hat{\beta}_T$ trained on the target samples available and another one $\M_{\T|\S}$ of estimated parameter $\hat{\beta}_{T|S}$ trained on the source samples and then enriched on the same target ones as $\M_\T$. We would like to know whether $\M_{\T|\S}$ will have a better prediction error than $\M_\T$ on a new sample $(x,y)$ drawn from the target distribution, i.e. if the transfer is positive or not. Let $f_{\hat{\beta}_T}(x)$ be the corresponding prediction from $\M_\T$ and $f_{\hat{\beta}_{T|S}}(x)$ be the one from $\M_{\T|\S}$. Following \cite{bosq2008inference}, we use the quadratic prediction error (QPE) as metric of the quality of a prediction $\RR (f_{\hat{\beta}}(x)) = \E[(y-f_{\hat{\beta}} (x))^2]$ where the expectation is taken with respect to the noise of $y$ and the distribution of $\hat{\beta}$. Therefore the transfer will be said negative for a given $x$ when $\E[(y-f_{\hat{\beta}_T} (x))^2] < \E[(y-f_{\hat{\beta}_{T|S}} (x))^2]$. Hence we would like to know in advance whether this inequality stands or not, i.e. when negative transfer will happen. In our work we derive a new quantity referred as \emph{gain} quantifying the benefits of transfer, without any assumption outside of the one of the linear model. While the hypothesis of a linear model may seem restrictive, it includes many variants such as generalized linear models (GAM) \cite{wood2017generalized} that make it possible to capture highly nonlinear effects through the use of spline bases. We will also show that it is possible to derive a hypothesis test to \emph{predict in practice} whether the transfer is positive or not. The contributions of the paper are the following:
\begin{enumerate}
    \item We formalize the problem of negative transfer for the fine-tuning of a linear regression model. However our framework is valid for a broad class of transfer procedures for the linear model found in the literature.
    \item We show that the transfer gain for a new feature vector $x$ depends on its representation on an eigenbasis depending on the parameters of the linear model.
    \item We establish a link between transfer by data pooling and fine-tuning and show that in the framework of linear regression they both yield estimators of the same form.
    \item 
    We suggest a statistical test to choose for a new observation $x$ between the target model $\M_\T$ or a fine-tuned one $\M_{\T|\S}$.
\end{enumerate}

The rest of the paper is organized as follows. Section 2 introduces the fine-tuning transfer procedure considered for the linear model and establishes the equation of the transfer gain. Section 3 presents the statistical test to predict positive and negative transfer. In Section 4 we apply the test on synthetic data as well as a real-world electricity consumption dataset. Section 5 illustrates how the gain relates to the source and target sample sizes. Finally Section 6 concludes our work and suggests further research possibilities. The appendices contain mathematical proofs and additional figures.
 
\section{Model transfer}
\label{theory}

Let $\D_S = \{(x_{S,i},y_{S,i}), i=1..N_S\}$ the source training dataset of size $N_S$ be, where the input vectors $x_{S,i} \in \R^D$ are supposed to be deterministic and the $y_{S,i}$ are supposed to be independent. We make the assumption of the gaussian linear model, i.e. that there exists $\beta_S \in \R^D$ such that $y_{S,i} = x_{S,i}^\top \beta_S + \veps_{S,i}$ with $\veps_{S,i} \sim \NN(0,\sigma_S^2)$. Similarly we define $\D_T = \{(x_{T,i},y_{T,i}), i=1..N_T\}$ the target training dataset with $x_{T,i} \in \R^D$ and independent $y_{T,i} = x_{T,i}^\top \beta_T + \veps_{T,i}$ where $\beta_T \in \R^D$ and $\veps_{T,i} \sim \NN(0,\sigma_T^2)$. The source and target data are supposed to be independent. The data can be rewritten more conveniently under matrix form $Y_\nu = X_\nu \beta_\nu + \veps_\nu$ where $\nu \in \{S,T\}$ denotes either the source or the target. The matrices $X_S$ and $X_T$ are the design matrices respectively of size $N_S \times D$ and $N_T \times D$ and are assumed to be full rank, so that $X_\nu^\top X_\nu$ is invertible. Hence it implies that $D<N_\nu$,which corresponds to low-dimensional setting. Many aspects of our work can be generalized to the high-dimensional setting, but is out of the scope of this paper. The standard procedure to estimate the coefficients is by minimization of the least-squares error $J_\nu(\beta) := \frac{1}{2} \norm{Y_\nu - X_\nu \beta}^2$ on each training set. It yields the well-known solution $\hat{\beta}_\nu = \Sigma_\nu^{-1} X_\nu^\top Y_\nu$ where $\Sigma_\nu = X_\nu^\top X_\nu$. In our framework we suppose that the number of samples $N_T$ to have a decent estimator of $\beta_T$ is too low, and that leveraging information from $\S$ may improve performance. %%Citer ?

\subsection{Fine-tuning of $\hat{\beta}_S$}

Therefore we start from the estimator $\hat{\beta}_S$  and fine-tune it by batch gradient descent (GD) of stepsize $\alpha$ on $\D_\T$. The following result gives the expression of the fine-tuned estimator.

\begin{proposition}
At iteration $k \in \N$ the fine-tuned estimator of $\beta_T$ is:
\begin{equation}\label{GD}
\hat{\beta}_k = A^k \hat{\beta}_S + (I-A^k) \hat{\beta}_T
\end{equation}
where $A=I_D - \alpha \Sigma_T$ and $I_D$ is the identity matrix of size $D$.
\end{proposition}

 We will refer to this resulting model as $\M_{\T|\S}$. Therefore the fine-tuned estimator is a matrix combination of source and target estimators. In fact this observation can be taken further in the right vector basis to give more insight on this expression. Since $\Sigma_T$ is symmetric and real-valued, let $P$ be an orthogonal diagonalization basis matrix such that $\Sigma_T = P \Lambda P^\top$ with $\Lambda = \textrm{diag}(\lambda_i, \, i=1..D)$ the diagonal matrix of eigenvalues of $\Sigma_T$. Let $\tilde{\beta}_\nu$ denote the coordinate of $\hat{\beta}_\nu$ in $\Sigma_T$'s eigenbasis. Hence $\hat{\beta}_\nu = P \tilde{\beta}_\nu$.  Reusing equation (\ref{GD}) yields:
\begin{equation}\label{GD_eigenbasis}
    \tilde{\beta}_k = (I_D-\alpha\Lambda)^k \tilde{\beta}_S + (I_D - (I_D-\alpha\Lambda)^k) \tilde{\beta}_T
\end{equation}

which means that for every coordinate $i$ in this basis we have:
\begin{equation}\label{GD_eigenbasis_coord}
    \tilde{\beta}_k^{(i)} = (1-\alpha \lambda_i)^k \tilde{\beta}_S^{(i)} + \big( 1 - (1-\alpha \lambda_i)^k \big) \tilde{\beta}_T^{(i)}.
\end{equation}

Hence when $\alpha$ is small enough and in the right basis, each coordinate of the fine-tuned coefficient is a convex combination of the source and target coefficients, albeit with different weights depending on the eigenvalues $\lambda_i$. For small eigenvalues of $\Sigma_T$ the fine-tuning procedure will give a larger weight to the source whereas it is the opposite for larger ones.

Note that these expressions relate this transfer strategy to the ones introduced in \cite{chen2015data}, where they consider two types of transfer for the linear model. The first one is the pooling of source and target data, leading to estimators of the form $\hat\beta_{\lambda} = W_\lambda \hat\beta_S + (I_D - W_\lambda) \hat\beta_T$ where $W_{\lambda}$ is a matrix depending on the penalty parameter $\lambda > 0$. The second one is a simple convex combination $\hat\beta(\omega) = \omega \hat\beta_S + (1-\omega) \hat\beta_T$ for a constant weight $\omega \in [0,1]$. Hence transfer by fine-tuning is between those two approaches: in the right basis and for $\alpha$ small enough each coefficient is a convex combination of the source and target ones, albeit with different weights depending on the eigenvalue $\lambda_i$. This allows for more adaptability than for a constant $\omega$ as will be shown in simulations. It is interesting to note that in the end two popular transfer approaches, namely data pooling and fine-tuning yield estimators of the same class $\hat\beta(W) = W \hat\beta_S + (I_D-W) \hat\beta_T$ with specific forms of $W \in \R^{d\times d}$. In the case of data pooling the expresion of $W$ is more complex and it is generally not symmetric (see \cite{chen2015data}). To our knowledge such a strong relationship between the approaches has never been highlighted in literature before.

%Two comments can be made on these expressions. The first one is that the idea of weighting estimators according to characteristics of the sample target population was already suggested in \cite{gao2008knowledge} and spontaneously appears with fine-tuning (although the weights are not dynamic). Depending on the nature of the target more or less weight will be given in certain directions.

\subsection{Transfer gain}

The quality of a model will be evaluated for a new independent sample $(x,y)$ drawn from the underlying distribution of $\T$. We want to know if for this given $x$ the estimator $\hat{\beta}_k$ learned on $\S$ but fine-tuned on $\T$ is better than the basic estimator $\hat{\beta}_T$.  Following what was discussed in the introduction we introduce the \emph{algebraic gain} $\Delta \RR_k (x)$ for sample $(x,y)$ defined by: $\Delta \RR_k (x) = \E[(y-\hat{y}_T)^2] - \E[(y-\hat{y}_k)^2]$ where $\hat{y}_T = x^\top \hat{\beta}_T$ and $\hat{y}_k = x^\top \hat{\beta}_k$. We have the following result in the case of fine-tuning.

\begin{proposition}
For transfer by fine-tuning as presented by equation (\ref{GD}), at iteration $k$ the gain is: \begin{equation}\label{gain_k}
    \Delta \RR_k(x)  =  x^\top H_k x \quad \textrm{where} \\
    H_k = \sigma_T^2 (\Sigma_T^{-1} - \alpha^2 \Omega_k \Sigma_T \Omega_k) - \sigma_S^2 A^k \Sigma_S^{-1} A^k  - A^k B A^k
\end{equation}
with $\Omega_k = \frac{1}{\alpha}\Sigma_T^{-1} (I_D - A^k)$, $B = (\beta_T - \beta_S)(\beta_T - \beta_S)^\top$. When it is positive, the transfer is beneficial for the sample $(x,y)$, and negative otherwise.
\end{proposition}

Therefore it can be seen that the matrix $H_{k}$ plays a significant role for the transfer problem. The gain will be positive for vectors in the span of the eigenvectors of $H_k$ associated to positive eigenvalues. The role of the noise in the data as well as the distance between the regression parameters also becomes clear with this formula and seems intuitive. When $\norm{\beta_S-\beta_T}$ is large, i.e. the means of $y_\nu$'s will differ significantly, transfer is likely to be detrimental. When $\sigma_T^2$ is large (the target data is noisy), the gain will increase since learning from the target data may be difficult. Note that this expression of the gain does not require any hypothesis on $x$, which is a major difference with previous works. We also see that an uniformly positive transfer may be impossible, and that the benefits of transfer are a local property: therefore for some $x$ it may be beneficial to use a fine-tuned model, whereas for others not. From (\ref{gain_k}) bounds on the prediction error can easily be derived:
\begin{equation}\label{oracle_ineq}
    \begin{array}{l}
         \E[(y-\hat{y}_k)^2] \le \E[(y-\hat{y}_T)^2] - \lambda_{\min} \,\,(H_k) \norm{x}^2 \\
         \\
         \E[(y-\hat{y}_k)^2] \ge \E[(y-\hat{y}_T)^2] - \lambda_{\max} \,(H_k) \norm{x}^2
    \end{array}
\end{equation}
where $\lambda_{\min}$ and $\lambda_{\max}$ respectively denote the minimum and maximum eigenvalues of $\Sigma_T$. Again those bounds do not require any assumptions and hold for any $x \in \R^D$ and only require $H_k$ to be symmetric, which is the case when performing transfer by fine-tuning. As one can see the transfer is always positive when $\lambda_{\min}(H_k) > 0$. More generally, a similar expression to (\ref{gain_k}) is possible for any estimator of the form $\hat{\beta}(W) = W \hat\beta_S + (I_D-W) \hat\beta_T$. However when $W$ is not symmetric, interpretability of transfer in terms of eigenvector direction is lost and the equations of (\ref{oracle_ineq}) cannot be established in the same way. Consequently if $H_k$ was accessible, one would know which model to use exactly. However the issue is that many quantities in the matrix are unknown, namely the true regression parameters $\beta_\nu$ and the true variances of the noise $\sigma_\nu^2$. A naive approach would consist to consider the "plug-in" estimate $\hat{H}_{k}$ by replacing the parameters by their estimates, but experiments have shown that this is a rather poor choice in most situations. Another strategy is therefore proposed in the next section. Finally we emphasize again that $x$ is potentially a \emph{novel} observation on which we require no hypothesis. In the aforementioned papers the bounds hold only under specific conditions that did now allow for any $x \in \R^D$, making our result broader.

\section{Statistical test for positiveness of transfer}
\label{theorie}

We simplify our problem to knowing in advance whether the transfer will be beneficial or not, i.e. if $\Delta \RR_k(x) >0$. Therefore an alternative is to define the problem as hypothesis testing. Considering that $\M_{\T|\S}$ is likely to be biased, we choose the null hypothesis $H_0 : \{ \Delta \RR_k(x) \le 0\}$ against the alternative $H_1 : \{ \Delta \RR_k(x) > 0\}$. This boils down to choosing between two models, the pure target one and the fine-tuned one for a given target sample. The idea of achieving best performance in transfer learning by taking advantage of multiple models could be related to \cite{gao2008knowledge} where they weighted classifiers according to the local properties of target observations.

\subsection{Expression of the test}

The main result of the paper is the following: \\

\begin{theorem}
\label{pvalue_approx}

Let $x \in \R^D$ be \textbf{any} observation. Let $\hat{\sigma}_S^2$ and $\hat{\sigma}_T^2$ be the estimations of the noise variances defined by $\hsigma_\nu^2 = \norm{Y_\nu - X_\nu \hat{\beta}_\nu}^2 / (N_\nu - D)$. Let $\rho$ be such that $\rho \ge \norm{\beta_T - \beta_S}/\sigma_T$. Then the following test is of approximate level $a$ to test $H_0$ against $H_1$:
\begin{equation}
 \mathbbm{1} \Big( \underbrace{\frac{\hsigma_T^2}{\hsigma_S^2} \frac{x^\top (\Sigma_T^{-1} - \alpha^2 \Omega_k \Sigma_T \Omega_k) x - \rho^2 \norm{A^k x}^2}{x^\top A^k \Sigma_S^{-1} A^k x}}_{:= \psi_k(x)} > q^{1-a} \Big)   
\end{equation}
where $q^{1-a}$ is the quantile of order $1-a$ of the $\F(N_T-D,N_S-D)$ Fisher-Snedecor distribution of degrees of freedom $N_T-D$ and $N_S-D$. The p-value for the observed data is:
\begin{equation}
    p_k(x) = \P_{F \sim \F(N_T-D,N_S-D)} \Big( F \ge \psi_k(x) \Big)
\end{equation}

\end{theorem}

The parameter $\rho$ can be seen as a prior on the distance between the source and target distributions. Indeed, in the gaussian case one can easily prove that $2 D_{KL}\Big(\NN(x^\top \beta_S,\sigma_S^2) || \NN(x^\top \beta_T,\sigma_T^2) \Big) \le g \big(\frac{\sigma_S^2}{\sigma_T^2} \big) + \rho^2 \norm{x}^2$ where $D_{KL}$ denotes the Kullback-Leiber (KL) divergence and $g(u) = u - \log(u) - 1$. The larger $\rho$ is, the more significant the difference between source and target distributions is allowed to be and thus the less likely the transfer will be beneficial. When $\rho = 0$ (i.e. $\beta_S = \beta_T$) only the variances differ. Note that the Cauchy-Schwarz approximation lowers the power of the test (see supplementary materials for more details). An issue is that $p_{k}(x) \longrightarrow 0$ when $k \rightarrow + \infty$. Hence when the number of gradient iterations goes to infinity, the test will almost systematically reject the null hypothesis, despite the gain converging to 0. Therefore a choice of a reasonable $k$ is of crucial importance. Finally the test can only be obtained when using a symmetric weight matrix $W$. \\

%The $\rho$ parameter can be seen as a prior on the distance $\norm{\beta_T - \beta_S}$ and the intensity of the target noise. The larger it is (i.e. the higher the discrepancy between source and target or the lower the target noise is), the higher the p-value and thus the lower the transfer is likely to be beneficial. 

%It is expected that for certain $x$ it is better to use $\hat{\beta}_T$, whereas in other cases $\hat{\beta}_k$ may be preferred. The user should therefore apply the test on the new feature vector: if the p-value $p_{k}(x)$ is low, $\M_{\T|\S}$ should be used for prediction. If it is high, then keeping $\M_\T$ is safer. This therefore defines a new model $\M_{\T|\S}^*$ equal to $\M_\T$ when $p_k(x) \ge 0.05$ (null hypothesis kept) and $\M_{\T|\S}$ when $p_k(x) < 0.05$ (null hypothesis rejected).

\subsection{Choice of $\alpha$, $k$ and $\rho$}
\label{tuning}

\begin{figure*}[h]
    \centering
    \subfloat[\label{fig:poly_overlap}][Overlap of $P_T$ and its different estimates.]{\includegraphics[width=0.44\linewidth]{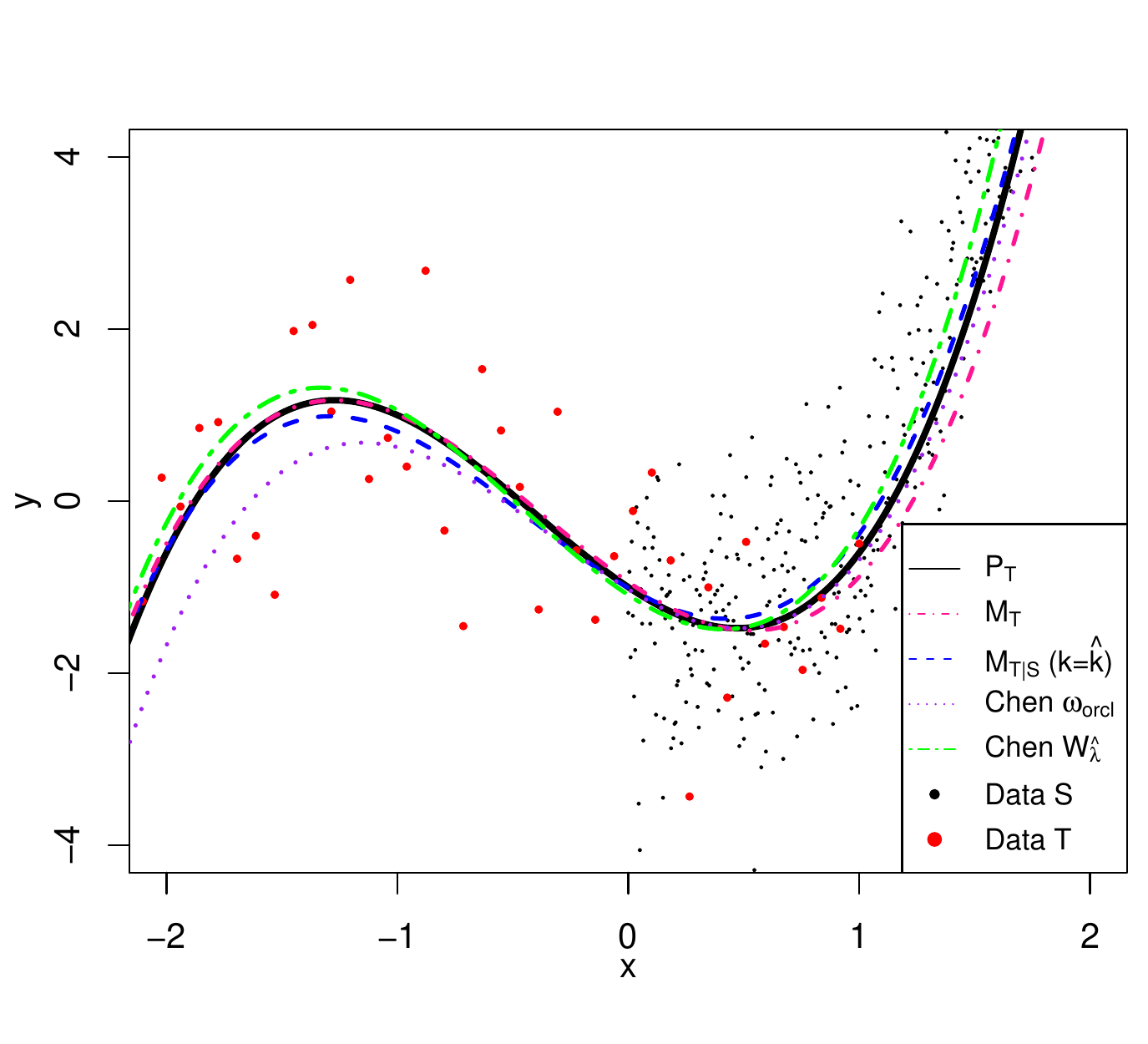}}
    \subfloat[\label{fig:poly_theo_gain}][Theoretical gain for different $k$]{\includegraphics[width=0.44\linewidth]{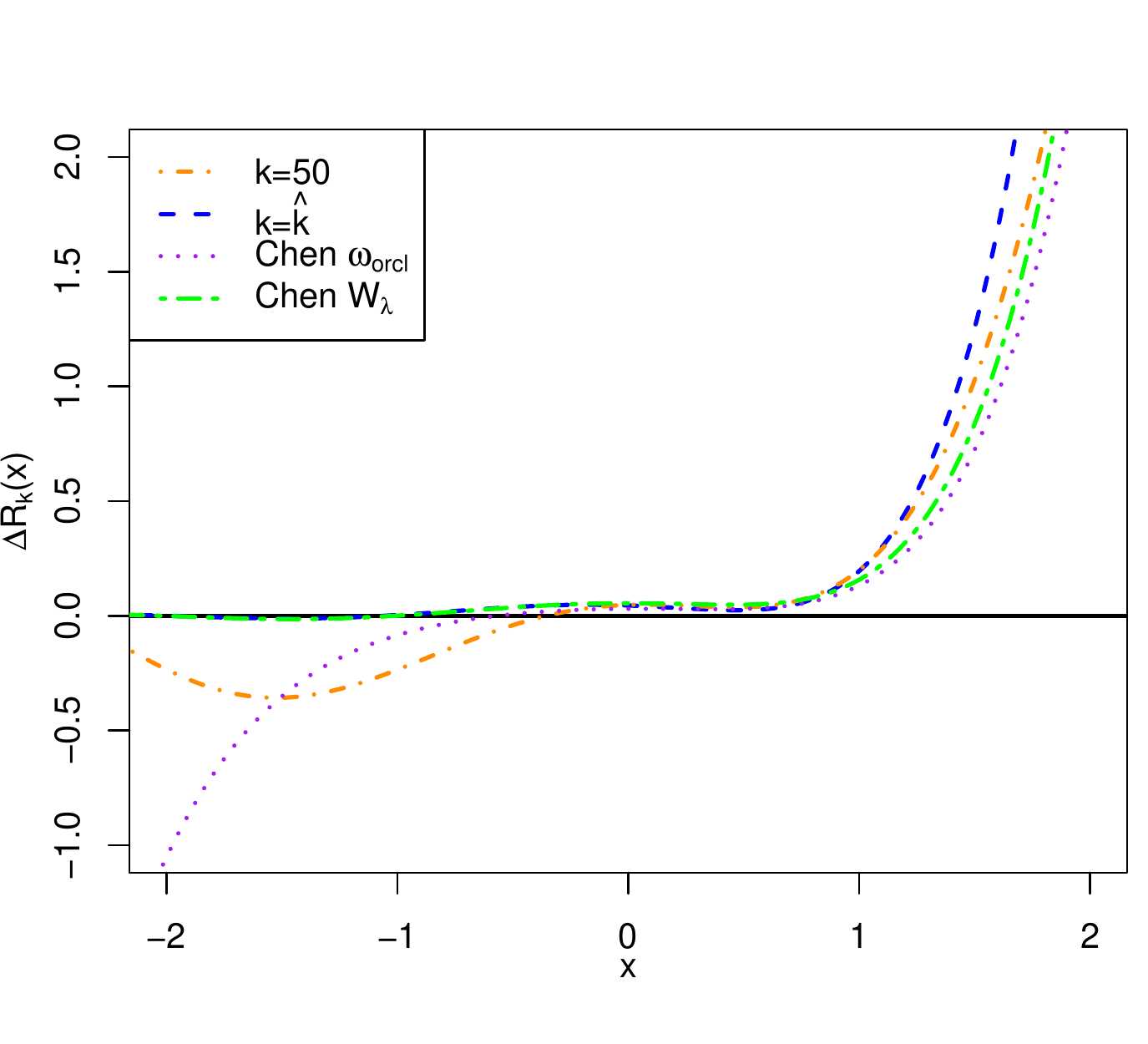}}
    \caption{Comparison of the estimates of $P_T$ and the theoretical gain values for the different models.}
    \label{fig:poly_overlap_estim_gain}
\end{figure*}

Three quantities must be tuned before usage of the test: the gradient step size $\alpha$, the number of iterations $k$ and the approximation parameter $\rho$. Equation (\ref{GD_eigenbasis_coord}) suggests that $0 < \alpha < 1 / \lambda_{\max}(\Sigma_T)$ so that the coordinates remain a convex combination of the source and target ones. Additionally according to \cite{bertsekas2015convex}, a step size $\alpha^* = 2/\big(\lambda_{\max}(\Sigma_T) + \lambda_{\min}(\Sigma_T) \big)$ allows to converge at optimal speed. However in our case convergence to $\hat{\beta}_T$ is not desirable since it would erase benefits from $\hat{\beta}_S$. Taking $\alpha = \alpha^* / 5$ or $\alpha^* / 10$ has proven to be a good choice in practice since it ensures the first condition while remaining close to $\alpha^*$. Moreover experimentally we observed that a low value of $\alpha$ could be compensated by a larger $k$. 

Ideally, one would choose the smallest $k$ such that $\lambda_{\min}(H_k) \ge 0$ (which would ensure an exclusively positive gain). However it depends on unknown parameters, and again a plug-in estimate yields poor results. The following approach yields satisfactory results. Let us denote by $\NN_T = \NN(x^\top \beta_T;\sigma_T^2 x^\top \Sigma_T^{-1} x)$ and $\NN_k =  \NN(x^\top \beta_{k} ; x^\top V_k x)$ the distributions of the predictions, where $\beta_{k} = \E[\hat{\beta}_k]$ and $V_k=\sigma_S^2 A^k \Sigma_S^{-1} A^k  + \sigma_T^2 \alpha^2 \Omega_k \Sigma_T \Omega_k$. It can be proved that:
\begin{equation}\label{eq:KL}
\Delta \RR_k (x) = -2\sigma_T^2 \, x^\top \Sigma_T^{-1} x \, D_{KL}(\NN_k || \NN_T) - \\ \sigma_T^2 \, x^\top \Sigma_T^{-1} x \ln \Big( \frac{x^\top V_k x}{\sigma_T^2 \, x^\top \Sigma_T^{-1} x} \Big)
\end{equation}
Let $U_k(x)$ denote the second term of the right hand side. Since the KL divergence is positive, it is needed that $U_k(x)$ is large to ensure a positive gain. Therefore maximizing $U_k(x)$ is more likely to maximize $\Delta \RR_k(x)$. Finally since the amount of target data is limited, we cannot afford to perform this procedure on a hold-out set. Therefore $k$ is selected by maximizing $\overline{U}_k := \frac{1}{N_S+N_T} \sum_{i=1}^{N_T+N_S} U_k(x_i)$ where the true variances have been replaced by their empirical counterparts. In case of absence of a local maximum, the elbow rule is applied instead. Finally, the choice of $\rho$ is performed by considering a range of possible values (typically between $10^{-5}$ and $1$) and checking the precision and recall of the test when used on the joint training data $\D_S \cup \D_T$. We refer by $\hat{k}$ and $\hat{\rho}$ the choices of $k$ and $\rho$ made with this procedure.

\section{Experiments}
\label{experiments}

\begin{figure*}[h]
    \centering
    \subfloat[][Test RMSE in function of $k$]{\includegraphics[width=0.44\linewidth]{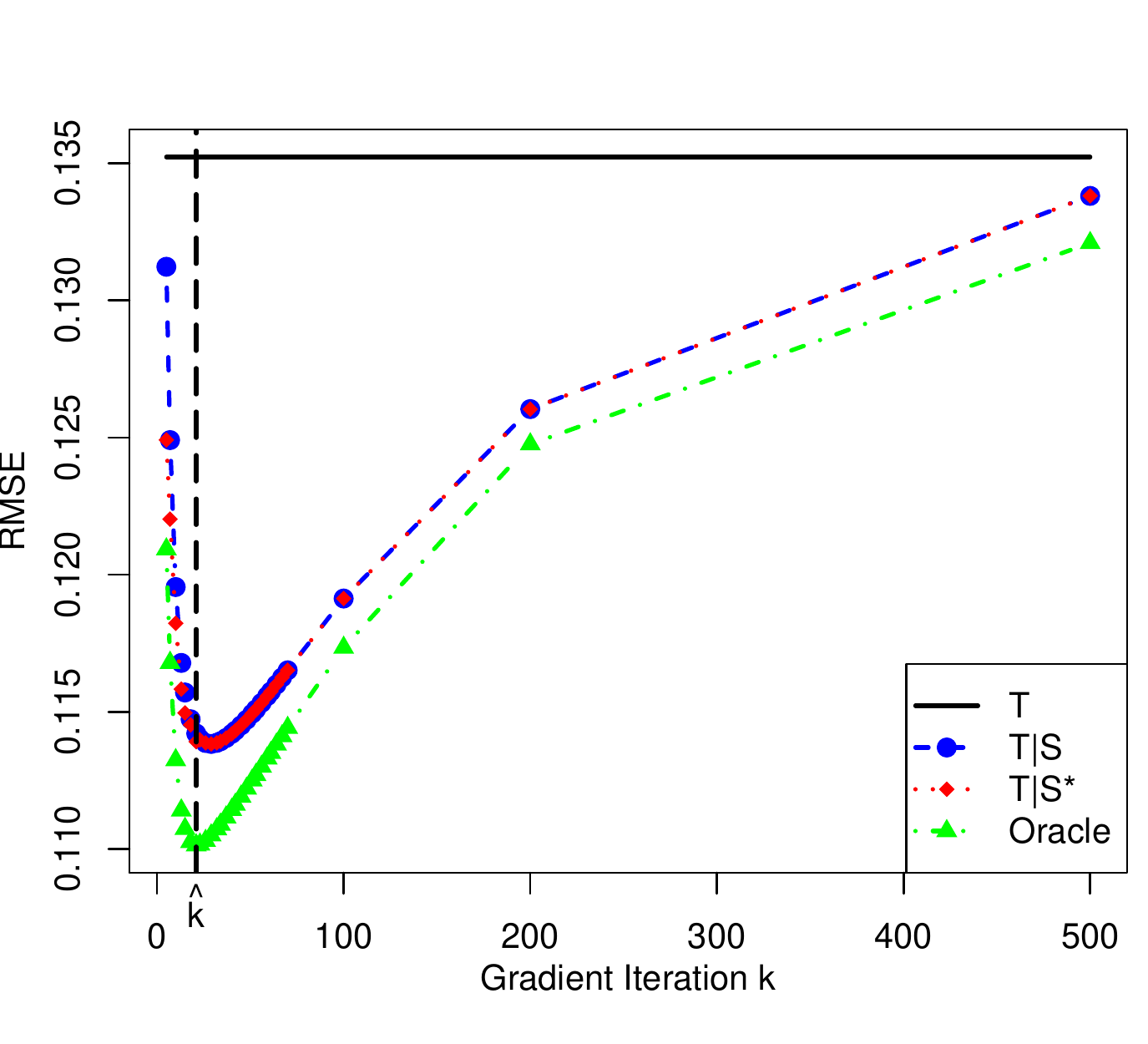}}
    \subfloat[][P-value over time on the test set ($k=\hat{k}$).]{\includegraphics[width=0.44\linewidth]{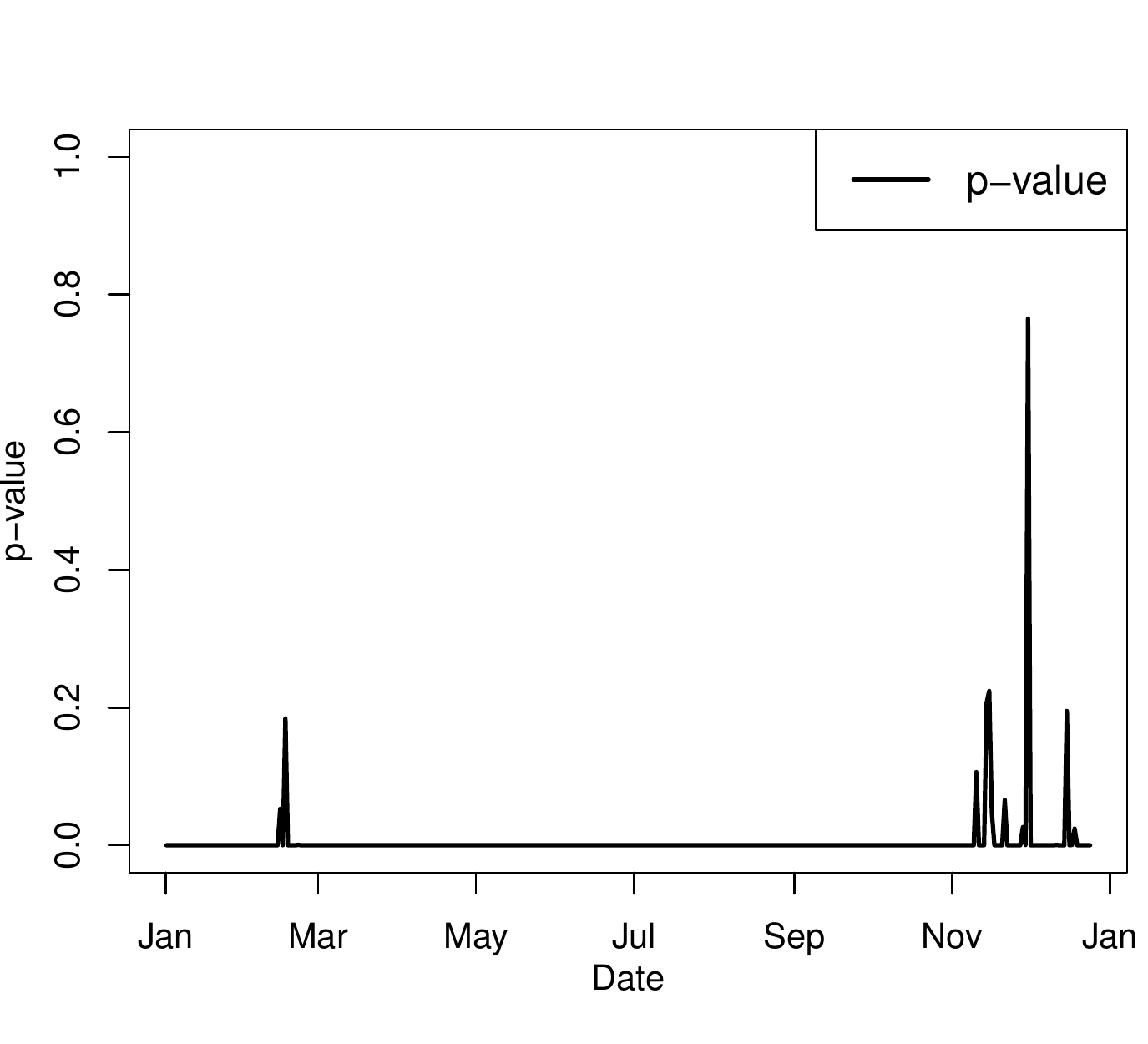}}
    \caption{Results on the test data (first scenario).}
    \label{fig:test_RMSE_A}
\end{figure*}

In this section the theoretical results are illustrated by numerical experiments. The first one is performed on synthetic data where all the parameters are known. In the second example real data of electricity consumption is used, proving the use of the test for real-life applications. %Finally we also conduct a simulation to show the dependence of the gain with the sample sizes $N_S$ and $N_T$.

\subsection{Low dimensional simulated data}

\begin{figure}[h]
    \centering
    \includegraphics[scale=0.44]{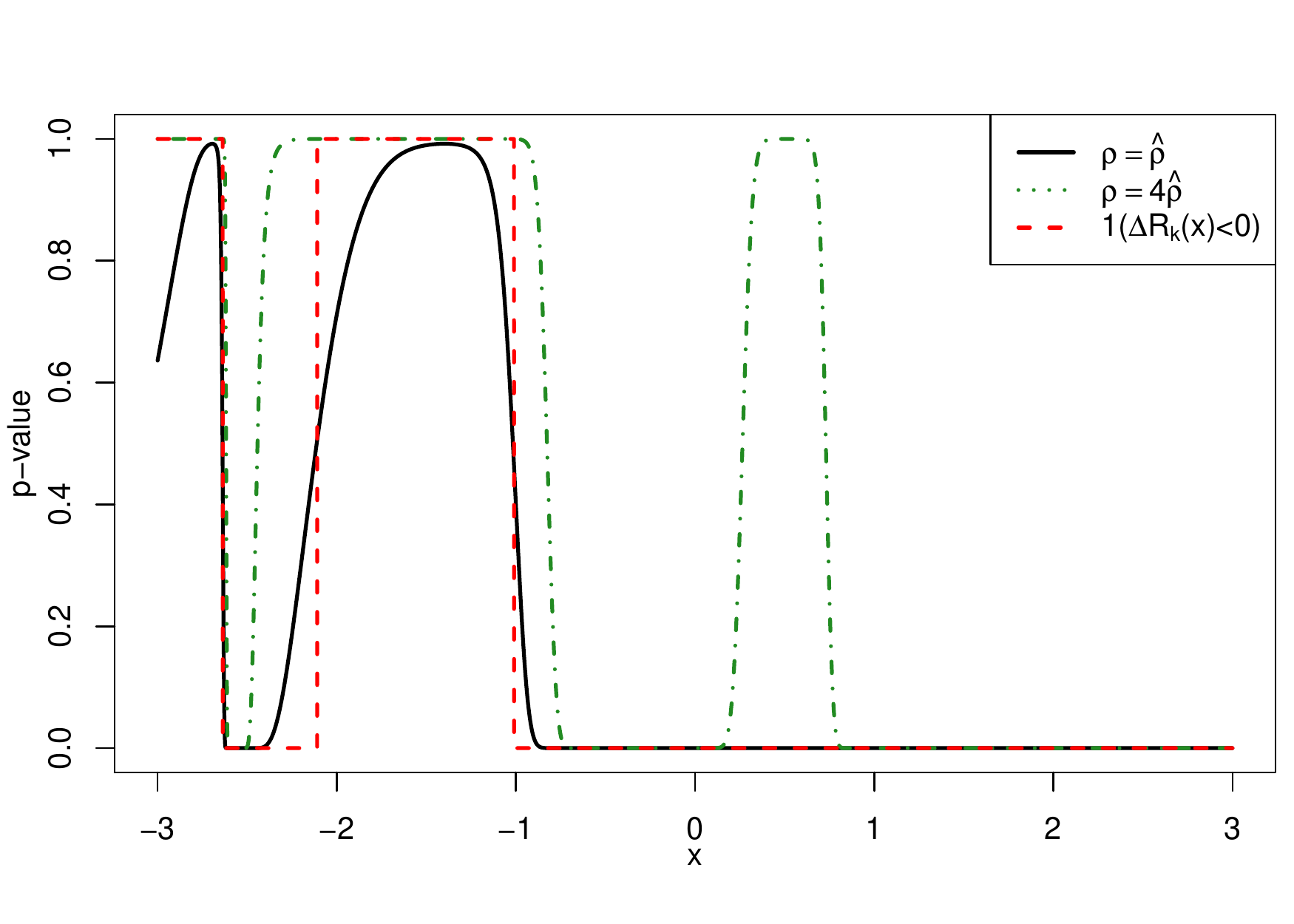}
    \caption{P-value in function of $x$ ($k = \hat{k}$).}
    \label{fig:poly_pval}
\end{figure}

We consider the problem of the estimation of the coefficients of a target polynomial $P_T(x) = \beta_{T,0} + \beta_{T,1} x + \beta_{T,2} x^2 + \beta_{T,3} x^3$ where $\beta_T = (-1,-1.8,1.2,1)^\top$. The advantage of this example lies in how it can be visualized, as one will see afterwards. We have $N_T = 60$ independent target observations $y_{T,i} = P_T(x_{T,i}) + \varepsilon_{T,i}$ with $x_{T,i} \in [-3,1]$ and $\varepsilon_{T,i} \sim \NN(0,\sigma_T^2)$. Additionally we have $N_S = 600$ independent source observations $y_{S,i} = P_S(x_{S,i}) + \varepsilon_{S,i}$ with $x_{S,i} \in [0,3]$ and $\varepsilon_{S,i} \sim \NN(0,\sigma_S^2)$. The coefficients of $P_S$ are the ones of $P_T$ plus a gaussian noise of mean $0$ and standard-deviation $0.3$, and set $\sigma_T^2 = \sigma_S^2 = 1$. Considering the locations of the samples for the source and target, the transfer is expected to be beneficial for $x \ge 1$. As suggested earlier the step size is set to $\alpha = \alpha^* / 10$, whereas the number of gradient iterations $k$ is decided by the strategy proposed in Section \ref{tuning}, leading to $\hat{k}= 405$. In order to illustrate the benefits of our tuning procedure, the results for $\hat{k}$ are compared with the ones of a sub-optimal $k=50$, but also the estimators using $W = \omega_{\textrm{orcl}} I_D$ or $W_{\hat\lambda}$ from \cite{chen2015data}. The true polynomial $P_T$ as well as its estimates are represented fig. \ref{fig:poly_overlap_estim_gain} (a). The fine-tuned model with $\hat{k}$ is the closest to the real curve, showing the advantages it takes from both the data of $\S$ and $\T$. Ergo those estimates confirm the aspect of the theoretical gain represented Figure \ref{fig:poly_overlap_estim_gain} (b). Furthermore fine-tuning proves to be superior to both of Chen's estimator, even with a theoretically optimal coefficient inaccessible in practice. This means that the fine-tuning procedure allows more flexibility in the use of both source and target data than simple data pooling or a convex combination of estimators. The benefits of the tuning procedure for $k$ is illustrated, as for $k=50$ the gain is still significantly negative between $-2$ and $-1$ whereas it is positive almost everywhere for $k=\hat{k}$. The p-values from (\ref{pvalue_approx}) are represented Figure \ref{fig:poly_pval} in function of $x$ for two different values of $\rho$. Additionally we represented the indicator $\mathbbm{1}\Big( \Delta \RR_k(x) \le 0\Big)$, i.e. when the theoretical gain is \emph{negative}. As one can see, the ranges of large p-values (typically greater than $0.05$) and the indicator usually coincide. The choice of $\rho = \hat{\rho}$ yields the best results, since for $\rho = 4\hat{\rho}$ the p-value suddenly increases between 0 and 1 despite the gain being positive.

\subsection{GEFCOM2012 electricity consumption}

\begin{figure}[h]
    \centering
    \subfloat[][Test RMSE in function of $k$.]{\includegraphics[width=0.44\linewidth]{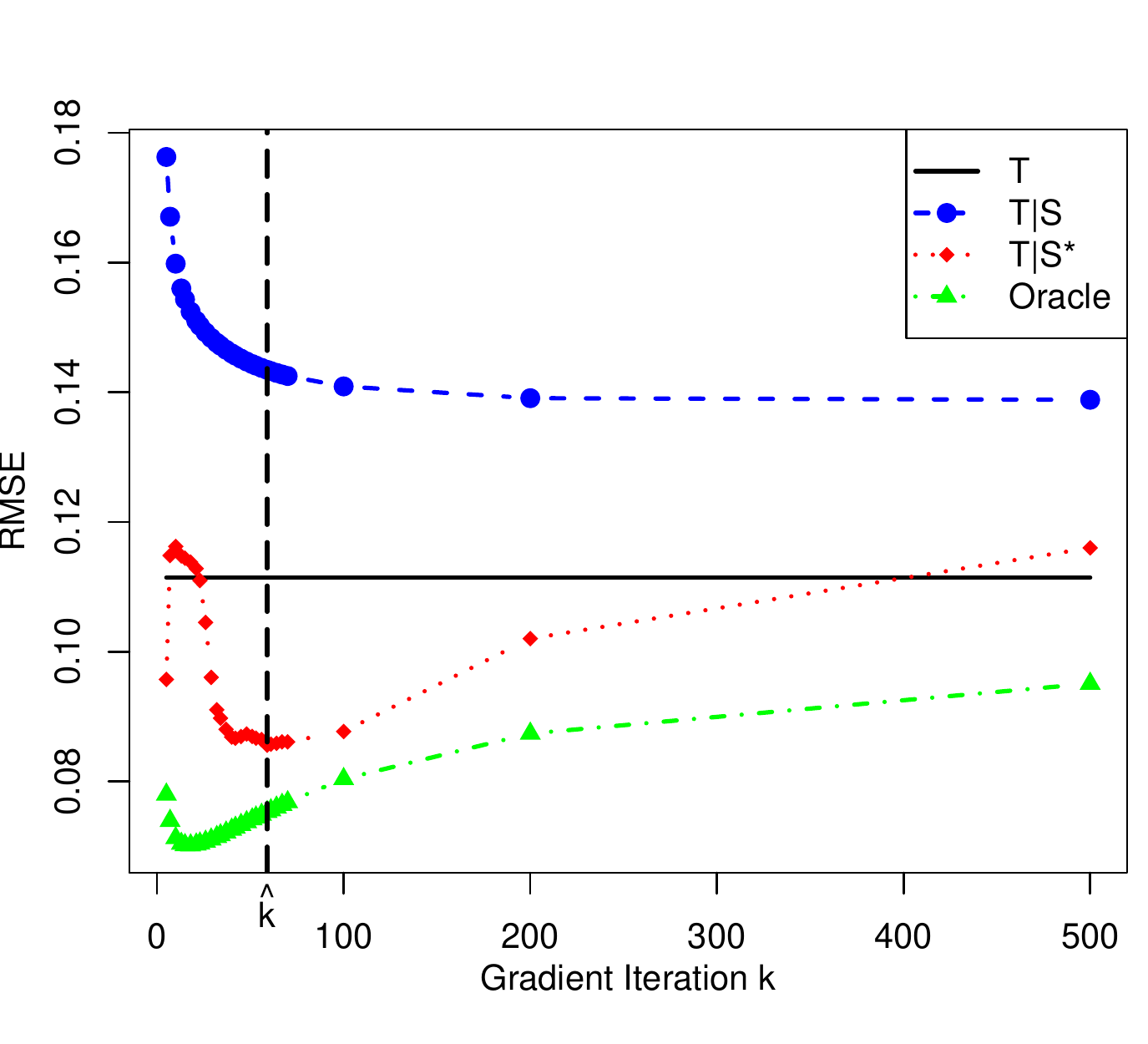}}
    \subfloat[][P-value over time on the test set ($k=\hat{k}$).]{\includegraphics[width=0.44\linewidth]{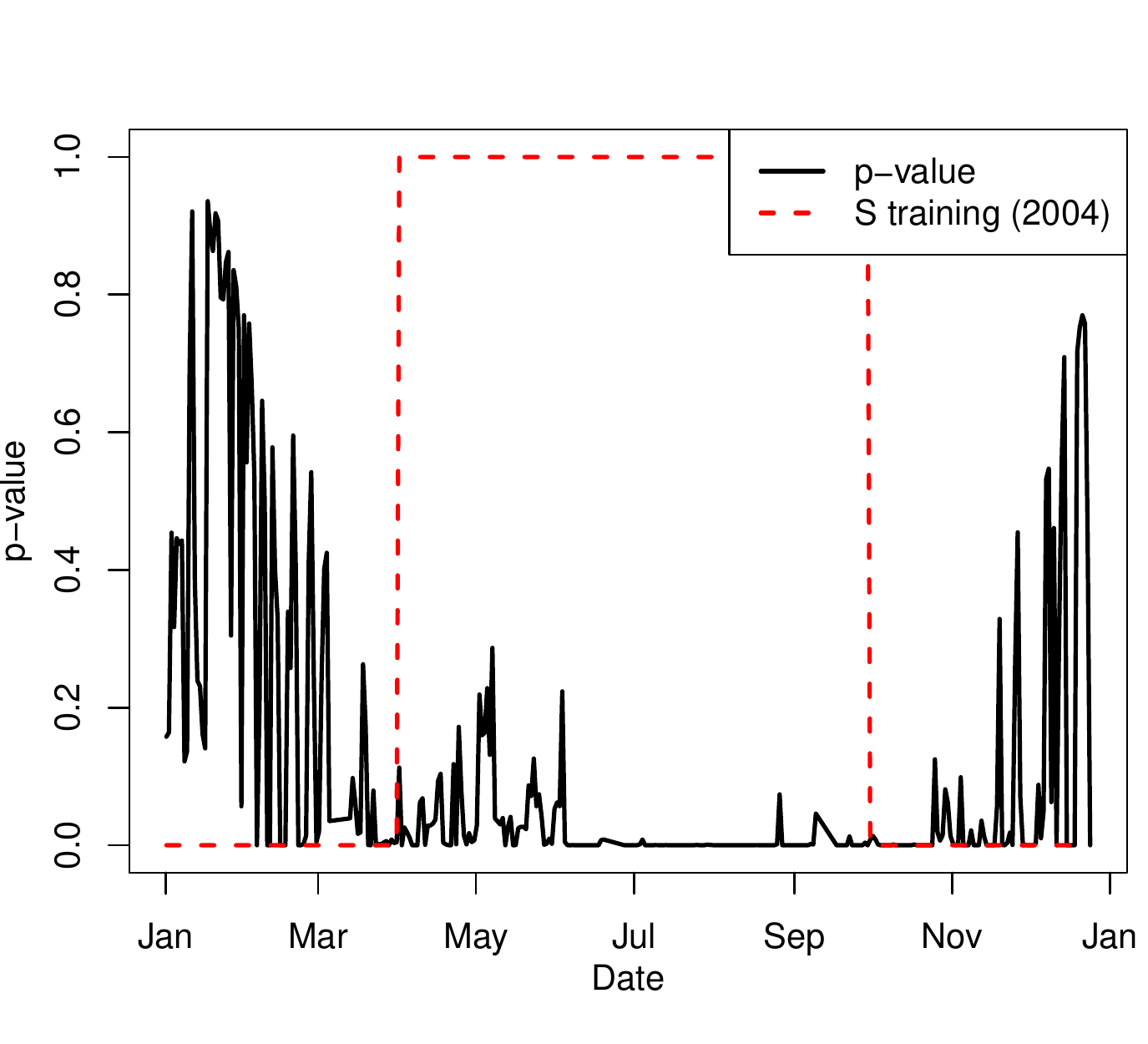}}
    \caption{Results on the test data (second scenario).}
    \label{fig:test_RMSE_B}
\end{figure}

This dataset was used during the GEFCOM2012 electricity consumption forecasting competition \cite{hong2014global}. It consists of the electricity consumption of 21 areas (zones) located in the United-States available from the 1\textsuperscript{st} of January 2004 to the 31\textsuperscript{st} of December 2007 with a 1 hour temporal resolution that we normalized. Input variables include calendar ones such as the day of the week, time of the year, but also the temperature measurements of 10 meteorological stations over the same period. The nature of our transfer is twofold: across time (period on which a model has been trained) and space (from one area to another). We will use the measured load at 8a.m. of zone 13 as source $\S$ and zone 2 as target $\T$. To focus on the benefits of our test and to avoid time series stationarity issues, both load time series have been detrended. The corresponding time series are represented in the supplementary material. In our work we focus on the use of our hypothesis test and not achieving pure predictive performance. Therefore the model we consider for both source and target is very simple:
\begin{equation}\label{model}
    y_{\nu,t} = \beta_{\nu,0} + \beta_{\nu,1} |\sin(\omega t)| + \beta_{\nu,2}  WE_t + \\ \sum_{j=3}^5 \beta_{\nu,j} \theta_t \, \mathbbm{1}\big( \theta_t \in I_j \big) + \varepsilon_{\nu,t}, \quad \nu \in \{S,T\}
\end{equation}
where $y_{\nu,t}$ is the load demand at 8a.m. for day $t$, $\omega = \frac{2\pi}{365}$. The sine term is used for the annual periodicity, $WE_t$ is a binary variable whose value is on $1$ on weekends. $\theta_t$ is the temperature and its effect has been cut into three intervals to translate the impact of heating and cooling on the electricity demand \cite{pierrot2011short}.  Whether it is the source or the target data, the training data will be included within the year 2004, whereas the test data on which performance is finally evaluated will be the whole 2005 year of zone 2.  \\
In order to evaluate the performance brought by our test, we consider the oracle prediction which knows in advance whether to use $\M_\T$ or $\M_{\T|\S}$. Thus the closer $\M_{\T|\S}^*$ is to the oracle, the better. The metric of evaluation will be the root mean squared error (RMSE) defined as: $RMSE = \sqrt{\frac{1}{T} \sum_{t=1}^T (y_t - \hat{y}_t)^2}$ where $T$ the number of test samples.

\subsubsection{First scenario}

In this scenario we suppose that the data for the source $\S$ is available for the whole year 2004. The target training data will only be available from October the 1\textsuperscript{st} to the end of the year. Hence $N_S=366$ and $N_T=92$. The RMSE for different values of $k$ is represented fig. \ref{fig:test_RMSE_A} (a), with a vertical line corresponding to our chosen $\hat{k}$. Here the improvement brought by the test is only marginal, for $k$ below a threshold. This is due to the discussed phenomenon at the end of Section \ref{theorie} where when $k \rightarrow \infty$ the test tends to systematically reject $H_0$. Note that the RMSE is minimal for $\hat{k}$. The errors for $\hat\omega_{\textrm{plug}}$ and $W_{\hat\lambda}$ from Chen et al. have been calculated but not represented because both always fared poorer than $\M_{\T}$. However most importantly the model $\M_{\T|\S}^*$ is always as  good as $\M_{\T|\S}$: the test can thus be used safely in practice. The p-value over time on the test set is also represented fig. \ref{fig:test_RMSE_A} (b). One sees that it's almost always close to 0, except locally for cold months. Since $\M_\T$ was trained on a similar period the year before, such a behavior is logical.

\subsubsection{Second scenario}

We consider the case where the training data from $\S$ is available between April the 1\textsuperscript{st} and September the 30\textsuperscript{th} 2004. The training data from $\T$ is available between the 1\textsuperscript{st} of September and the 31\textsuperscript{st} of December 2004, and thus $N_S=182$ and $N_T=122$. In practice it could correspond to the case where a customer breaks his contract, and a new one arrives.  

Results on the test data for the different approaches are given fig. \ref{fig:test_RMSE_B} (a). We see that this time the test significantly improves upon the individual forecasts, lowering the RMSE by more than 0.02 compared to $\M_\T$ for the chosen $\hat{k}$. The test efficiently detects the situations of positive and negative transfer, thus taking advantage of each model's specificities. Furthermore the prediction for $\hat{k}$ is very close to the oracle. Again both of Chen et al.'s predictors yielded poor results, being only marginally better than $\M_{\T}$. The p-value over time on the test set is also plotted fig. \ref{fig:test_RMSE_B} (b). It is close to 0 on a period similar to the one the source model was trained the year before, and large during the cold months where the model $\T$ is expected to be better.

\section{Interpretation of the gain with sample sizes}

\begin{figure*}[h]
    \centering
    \subfloat[\label{fig:poly_overlap}][$k=0$]{\includegraphics[width=0.33\linewidth]{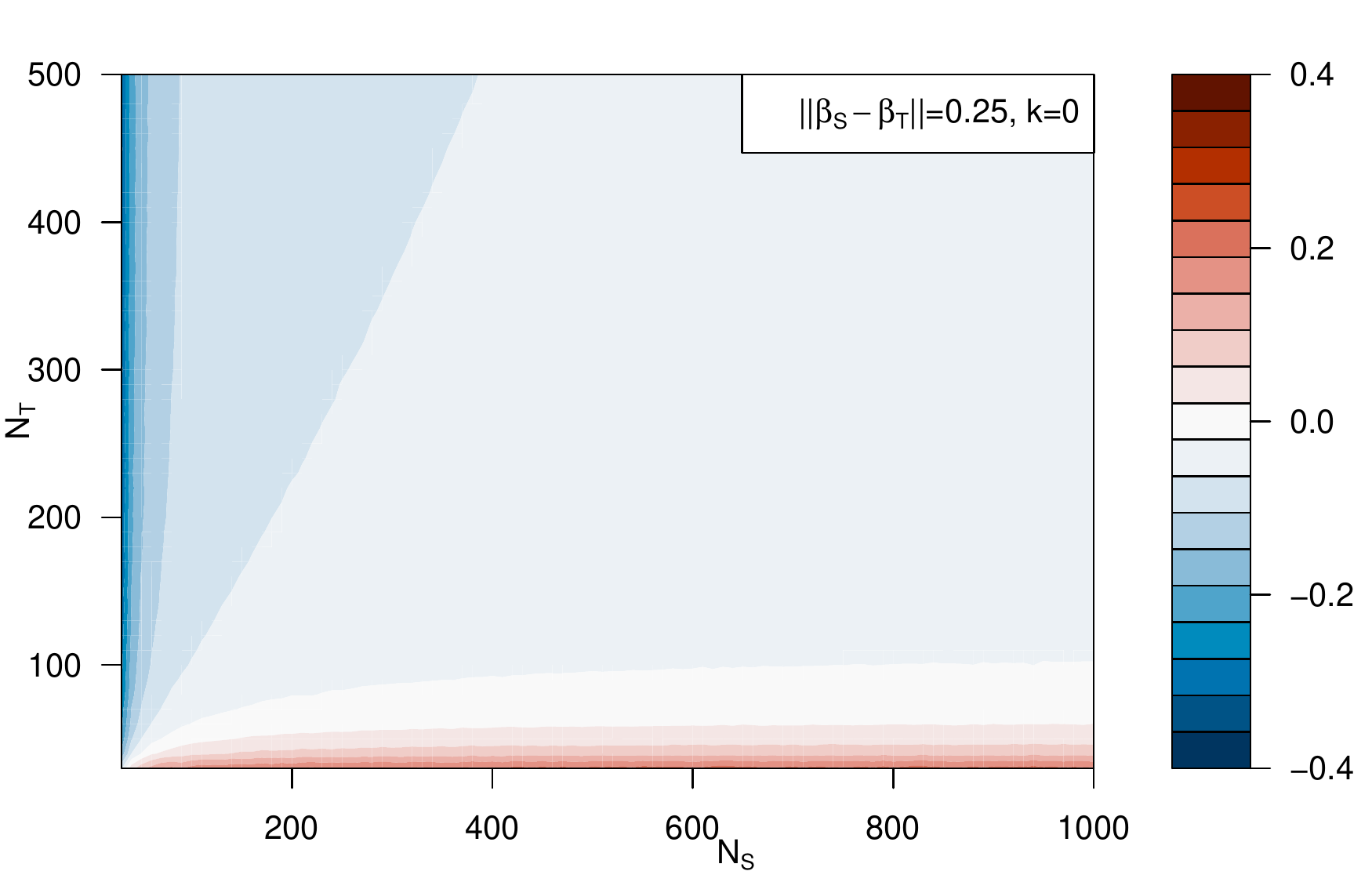}}
    \subfloat[\label{fig:poly_theo_gain}][$k=10$]{\includegraphics[width=0.33\linewidth]{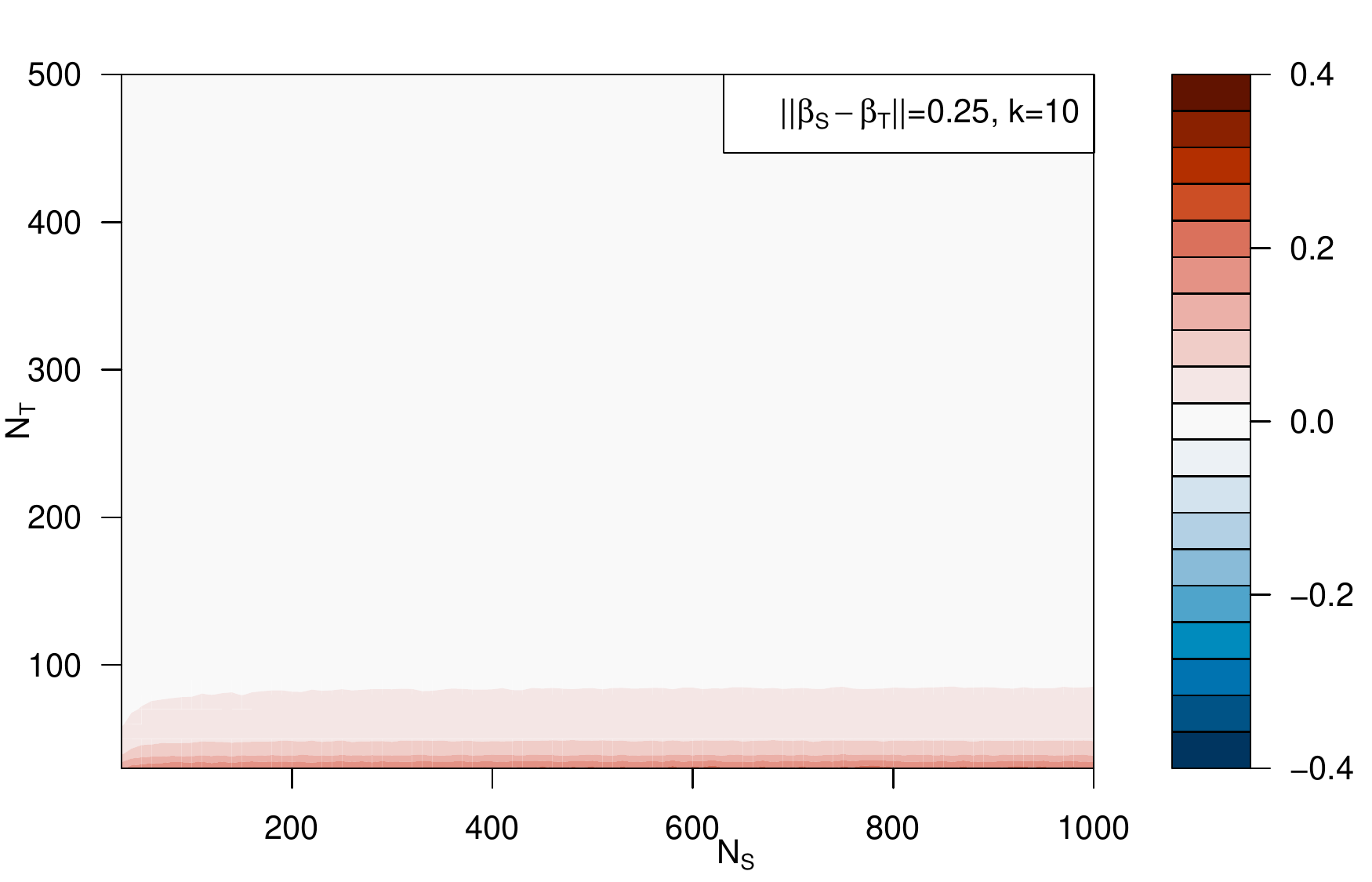}}
    \subfloat[\label{fig:poly_theo_gain}][$k=50$]{\includegraphics[width=0.33\linewidth]{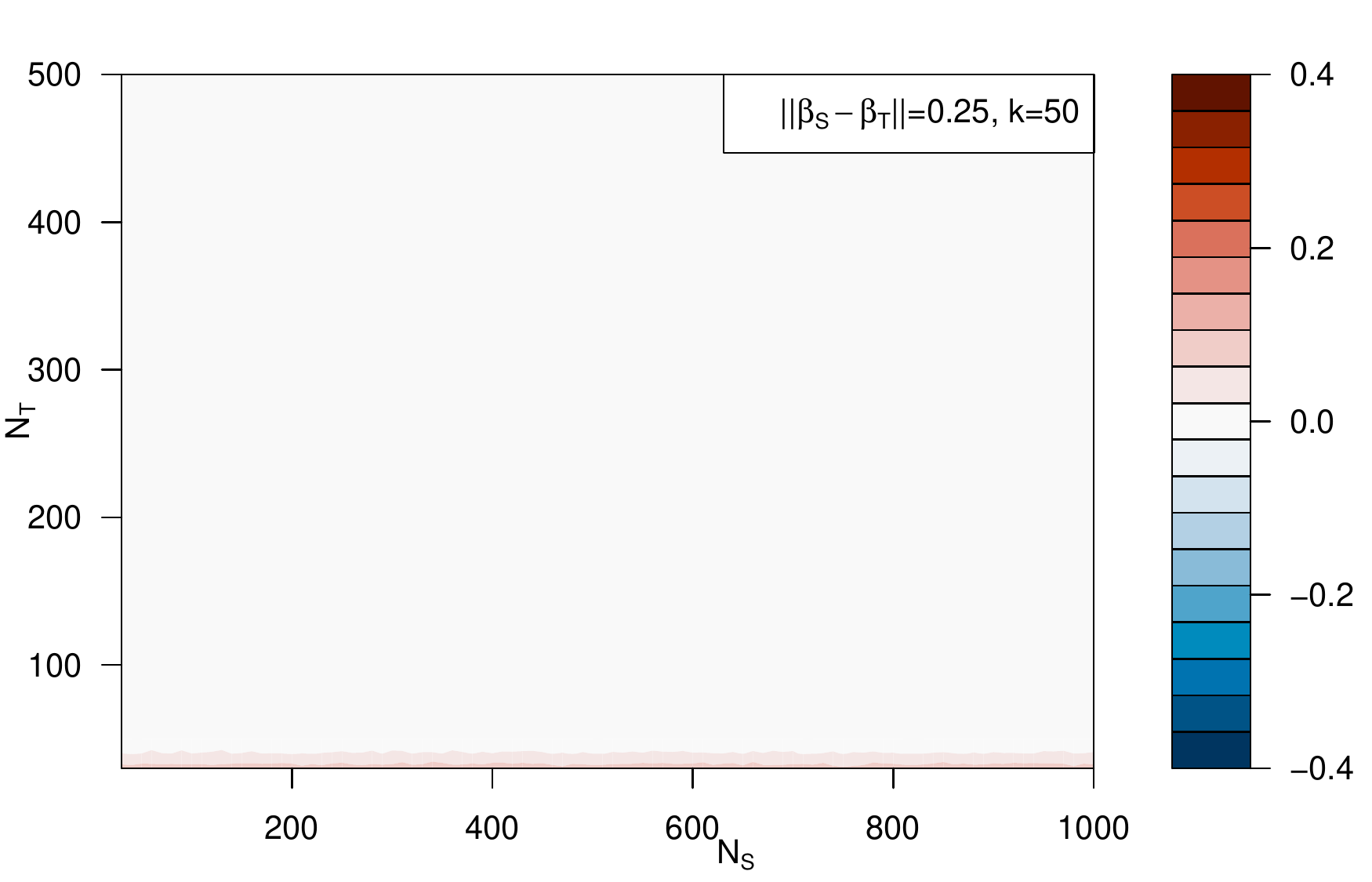}}
    \caption{Transfer phases in function of $N_S, N_T$ and $k$.}
    \label{fig:gain_phases}
\end{figure*}

A natural question to ask is how the gain evolves with the sample sizes $N_S$ and $N_T$. One would for instance expect the gain to increase when the number of source samples is order of magnitudes higher than the target one. In order to analyze these dependencies, we consider the following experimental framework. We suppose that the source and target data are i.i.d. $x_{\nu,i} \sim \NN (0,I_D)$ (thus $\Sigma_\nu \sim \mathcal{W}_D(I_D,N_\nu)$ where $\mathcal{W}_D(I_D,\Psi)$). For $(N_S,N_T)$ in a grid $I_S \times I_T$, we calculate and average the gain $\Delta \RR_k(x)$ over $B=50$ simulations for $x \sim \NN(0,I_D)$ as well. Algorithm \ref{alg:simulation} summarizes the procedure. This experiment is conducted for a dimension size $D=15$, $k\in\{0,10,50\}$, $\alpha = \alpha^*/5$ and $\norm{\beta_S-\beta_T} = 0.25$ (both coefficients have been randomly sampled). In order to improve the readability, the gain has been thresholded to the range $[-0.4,0.4]$. The results are represented in Fig. \ref{fig:gain_phases}.\\

Phases are observed depending on the values of $N_S$ and $N_T$, and follow the intuition. When no fine-tuning is performed (i.e. $k=0$) the gain will be positive only when the number of target samples $N_T$ is small and the number of source ones $N_S$ is large enough. For $N_T$ above a certain threshold, negative transfer will systematically happen. When $k$ increases, the blue areas corresponding to negative gain fade away thus meaning that at worst the transfer procedure will have a neutral impact, even for large values of $N_T$. The benefits of transfer through fine-tuning are particularly visible for $k=10$ with an increase of the size of the positive transfer areas: the fine-tuning procedure allows to take advantage of both source and target samples. However as emphasized before, an excessive number of gradient iterations may erase the benefits of transfer as seen in Figure \ref{fig:gain_phases} (c) obtained for $k=50$ where only for extremely small values of $N_T$ transfer can be beneficial. This is because the fine-tuned estimator $\hat\beta_k$ has come too close to the pure target one $\hat\beta_T$. Note that these figures remind of Fig. 1 from \cite{ben2010theory}.

%If $k$ isn't too large (e.g. $k=10$), for low values of $N_T$ we still observe strong red areas: the fine-tuning procedure allows to take advantage of both source and target samples. However when the number of gradient iterations is too high, the red areas fade away. This comes because excessive gradient iterations erase the benefits of the estimator learned on the source data and $\hat{\beta}_k$ has become too close to $\Hat{\beta}_T$. Note that these figures remind of Fig. 1 from \cite{ben2010theory}

\begin{algorithm}[h]
%label
{\caption{Gain simulation for varying $N_S$ \& $N_T$}
\label{alg:simulation}}% caption
{% contents
{\bf Initialisation}: $D, \beta_\nu, \sigma_\nu^2,k$. $I_S =  \{30,40,\dots,1000\}$ and $I_T =  \{30,40,\dots,500\}$.

\medskip

{\bf Recursion}: For $(N_S,N_T)$ in $I_S \times I_T$:
\begin{enumerate}

\item $\overline{\Delta \RR_k}(N_S,N_T) \leftarrow 0$

\medskip

{\bf{Recursion} For $b=1,\dots,B$:}
\begin{enumerate}
    \item  Generate $X_\nu \sim \NN(0,I_D)$. Deduce $\Sigma_\nu$.
    \item Generate $x\sim \NN(0,I_D)$
    \item Calculate\\ 
    $\overline{\Delta \RR_k}(N_S,N_T) \leftarrow \overline{\Delta \RR_k}(N_S,N_T)+ (1/B) \, x^\top H_k x$.
\end{enumerate}

\end{enumerate}
}
\end{algorithm}

\section{Conclusion and future work}

In this paper a novel framework for the problem of transfer learning for the linear model is proposed. By defining the gain of transfer by a difference of quadratic prediction errors, we obtain a quantity that measures how beneficial or detrimental transfer by gradient descent is for a new (potentially unobserved) $x$. However the framework of the gain is applicable for any estimator of the form $\hat\beta(W) = W \hat\beta_S + (I_D-W) \hat\beta_T$, which encompasses many found in the literature. Since this gain depends on unknown parameters in practice, we derived a statistical test relying of the Fisher-Snedecor $\mathcal{F}$ distribution to predict negative transfer. The test was applied on synthetic as well as real-world electricity demand data, where it proved its ability to predict negative transfer for new observations. \\
However despite its success, some points remain to investigate. How to choose the right number of gradient iterations $k$ remains problematic, although an empirical approach has been suggested. Furthermore in order to obtain a tractable calculation and satisfying empirical results, we had to rely on an approximation. Another possibility would be to transfer only a subset of parameters. This is often the case for neural networks where only certain layers are transferred \cite{laptev2018reconstruction}, but could be adapted for linear models. \cite{dar2020double} investigate the benefits of transfer depending on the number of parameters transferred, but do not indicate how to choose the subset to transfer. Moreover they do it in the static setting, i.e. no fine-tuning. We have also supposed that the matrices $\Sigma_\nu$ are invertible. However defining the gain without this hypothesis is still possible although its form is slightly more complex, which makes it difficult to adapt the test directly. Finally in this paper we made the hypothesis of linearity, which could seem restrictive. However nonlinearity can be achieved through generalized additive models (GAM) for instance. Since they boil down to a linear model, the formula of the gain is valid for it as well. However as such, the test we introduced cannot be used with GAM yet, and how to extrapolate it is currently under investigation.

\appendix

\section*{Appendices}

The first section of the appendix presents proofs of the results of the main article (sections 1 through 3). The second part consists in figures relating to the experimental Section 4.

\section{Proofs details}

This appendix presents proofs of the results in the paper.

\subsection{Proposition 1}

\begin{proof}
We proceed by mathematical induction. 

\begin{itemize}
    \item For $k=0$ the property is trivial. $\hat{\beta}_0 = \hat{\beta}_S = A^0 \hat{\beta}_S + (I_D - A^0) \hat{\beta}_T$. 
    \item Let $k \in \N$ be. We suppose the property true at rank $k$. We have 
    
    $$
    \hat{\beta}_{k+1} = \hat{\beta}_k - \alpha \nabla J_T(\hat{\beta}_k) = \hat{\beta}_k - \alpha \Sigma_T \hat{\beta}_k + \alpha X_T^\top Y_T
    $$
    
    By definition of $A = I_D - \alpha \Sigma_T$ and because $X_T^\top Y_T = \Sigma_T \hat{\beta}_T$ we obtain:
    
    $$
    \hat{\beta}_{k+1} = A \hat{\beta}_k + \alpha \Sigma_T \hat{\beta}_T
    $$
    
    Finally by induction hypothesis:
    
    $$
    \hat{\beta}_{k+1} = A [A^k \hat{\beta}_S + (I_D - A^k) \hat{\beta}_T] + \alpha \Sigma_T \hat{\beta}_T = A^{k+1} \hat{\beta}_S + (I_D - A^{k+1})\hat{\beta}_T
    $$
    
\end{itemize}

which concludes the induction. 
\end{proof} %\qed

\subsection{Equations (\ref{GD_eigenbasis}) and (\ref{GD_eigenbasis_coord})}

\begin{proof}
Let $P$ be the orthogonal matrix of eigenvectors of $\Sigma_T$ be, i.e. such that $\Sigma_T = P \Lambda P^\top$ with $\Lambda = \textrm{diag}(\lambda_i, i=1..D)$ and $P P^\top = P^\top P = I_D$. Thus $\hat\beta_\nu = P \tilde\beta_\nu \Leftrightarrow \tilde\beta_\nu = P^\top \hat\beta_\nu$. One can also write that $A = P (I_D - \alpha\Lambda) P^\top$. Hence reinjecting in (1) gives:

\[
\hat\beta_k = P (I_D - \alpha\Lambda)^k P^\top \hat\beta_S + P (I_D - (I_D-\alpha\Lambda)^k) P^\top \hat\beta_T
\]

Applying $P^\top$ on the left of this equation yields:

\[
\tilde\beta_k = (I_D-\alpha\Lambda)^k \tilde\beta_S + (I_D - (I_D-\alpha \Lambda)^k) \tilde\beta_T
\]

Finally the matrices involved are diagonal with respective terms $(1-\alpha\lambda_i)^k$ and $1-(1-\alpha\lambda_i)^k$, thus resulting in equation (3).
\end{proof}

\subsection{Proof of Proposition 2}

\begin{proof}
We remind that $\hat{\beta}_S \sim \NN(\beta_S,\sigma_S^2 \Sigma_S^{-1})$ and $\hat{\beta}_T \sim \NN(\beta_T,\sigma_T^2 \Sigma_T^{-1} )$. By independence of $\hat{\beta}_S$ and $\hat{\beta}_T$ we thus have 

$$\hat{\beta}_k \sim \NN \Big(A^k \beta_S + (I_D - A^k) \beta_T,\sigma_S^2 A^k \Sigma_S^{-1} A^k + \sigma_T^2 (I_D - A^k) \Sigma_T^{-1} (I_D - A^k) \Big)$$

It is easy to see that $\sigma_T^2 (I_D - A^k) \Sigma_T^{-1} (I_D - A^k) = \sigma_T^2 \alpha^2 \Omega_k \Sigma_T \Omega_k$. We will note $\beta_k = \E[\hat{\beta}_k]$ and $V_k = \Var(\hat{\beta}_k)$. For an independent $y = x^\top \beta_T + \veps$ with $\veps \sim \NN(0,\sigma_T^2)$ we obtain that:

$$
\begin{array}{l}
y - \hat{y}_T = x^\top (\beta_T - \hat{\beta}_T) + \veps \sim \NN \Big(0, \sigma_T^2 (1+x^\top \Sigma_T^{-1} x) \Big) \\
\\
y-\hat{y}_k = x^\top (\beta_T - \hat{\beta}_k) + \veps \sim \NN \Big( x^\top (\beta_T - \beta_k) , \sigma_T^2 + x^\top V_k x\Big)
\end{array}
$$

Thus $\RR(\M_T) = \E[(y-\hat{y}_T)^2] = \Var(y - \hat{y}_T) + \E[y-\hat{y}_T]^2 = \sigma_T^2 (1+x^\top \Sigma_T^{-1} x)$ and $\RR(\M_{T|S}) = \E[(y-\hat{y}_k)^2] = \sigma_T^2 + x^\top V_k x + \Big(x^\top(\beta_T - \beta_k) \Big)^2
$. Therefore:

$$
\Delta \RR_k(x) = \sigma_T^2 x^\top \Sigma_T^{-1} x - x^\top V_k x - x^\top B_k x
$$

Finally noticing that $\beta_T - \beta_k = A^k(\beta_T - \beta_S)$, we obtain that $\Big(x^\top(\beta_T - \beta_k) \Big)^2 = x^\top A^k B A^k x$ where $B=(\beta_T - \beta_S)(\beta_T - \beta_S)^\top$ thus yielding the expected result. 
\end{proof}

\subsection{Proof of the equations of (\ref{oracle_ineq})}

\begin{proof}
$H_k$ is symmetric. Hence we can introduce $\{u_i\}_{i=1..D}$ an orthonormal basis of eigenvectors of it with $\lambda_i(H_k)$ the associated eigenvalues. Let $x \in \R^D$ be with coordinates $x_i$ in this basis. Thus $x$ can be rewritten $x = \sum_{i=1}^D x_i u_i$. Since $\{u_i\}$ is orthonormal (i.e. $u_i^\top u_j = 1$ if $i=j$ and $0$ else) it follows that:

\[
x^\top H_k x = \sum_{i,j=1}^D \lambda_i(H_k) \, x_i x_j \, u_i^\top u_j = \sum_{i=1}^D \lambda_i(H_k) \, x_i^2.
\]
 Since $\lambda_{\min}(H_k) \le \lambda_i(H_k) \le \lambda_{\max}(H_k)$ we get that $\lambda_{\min}(H_k) \norm{x}^2 \le x^\top H_k x \le \lambda_{\max} (H_k) \norm{x}^2$. Finally remembering that $x^\top H_k x = \E[(y-\hat{y}_T)^2]-\E[(y-\hat{y}_k)^2]$ yields (5).
\end{proof}

\subsection{Proof of Theorem 1}

\begin{proof}
It would be natural to reject $H_0$ if an estimator $\hat{\delta}(x)$ of the gain is above a certain threshold. Hence a natural form of such a decision rule is $\mathbbm{1}( \hat{\delta}(x) > K_a)$, where $K_a$ is a constant depending on the desired level $a$ of the test. We consider the estimator of $\Delta \RR_k(x)$:

\[
\hat{\delta}(x) = \hsigma_T^2 \, x^\top \big( \Sigma_T^{-1} - \alpha^2 \Omega_k \Sigma_T \Omega_k \big) x - \hsigma_S^2 x^\top A^k \Sigma_S^{-1} A^k x - x^\top A^k B A^k x
\]

While the matrix $B$ is not accessible in practice, we start from this estimator for the sake of the simplicity of the calculations. We will address this issue later. It can be proved (see hereafter) that the type I error, the probability of wrongly rejecting the null hypothesis, is the largest at the boundary $\Delta \RR_k(x) = 0$. Thus $\hat{\delta}(x) > K_a$ is equivalent to:

\[
\frac{\hsigma_T^2 / \sigma_T^2}{\hsigma_S^2 / \sigma_S^2}
 + \frac{\hsigma_T^2}{\hsigma_S^2} \, \frac{x^\top A^k B A^k x}{\sigma_T^2 x^\top A^k \Sigma_S^{-1} A^k x} > \frac{K_a + x^\top A^k B A^k x + \hsigma_S^2 x^\top A^k \Sigma_S^{-1} A^k x}{\hsigma_S^2 x^\top A^k \Sigma_S^{-1} A^k x}
\]

Since $\frac{\hsigma_T^2 / \sigma_T^2}{\hsigma_S^2 / \sigma_S^2} \sim \mathcal{F}(N_T-D,N_S-D)$, taking $K_a = q^{1-a} \hsigma_S^2 x^\top A^k \Sigma_S^{-1} A^k x - x^\top A^k B A^k x - \hsigma_S^2 x^\top A^k \Sigma_S^{-1} A^k x + \frac{\hsigma_T^2}{\hsigma_S^2} x^\top A^k B A^k x$  (where $q^{1-a}$ is the quantile of order $1-a$ of the $\mathcal{F}(N_T-D,N_S-D)$ distribution) yields the test of level $a$:

\[
\mathbbm{1} \Big( \phi_k(x) := \frac{\hsigma_T^2 x^\top \big(\Sigma_T^{-1} - \alpha^2 \Omega_k \Sigma_T \Omega_k\big) x - (\hsigma_T / \sigma_T)^2 x^\top A^k B A^k x}{\hsigma_S^2 x^\top  A^k \Sigma_S^{-1} A^k x} > q^{1-a} \Big)
\]

However $B$ and $\sigma_T$ are unknown in practice, we will thus have to rely on a lower bound of $\phi_k(x)$ for the test. By hypothesis, we have $\norm{\beta_T - \beta_S} / \sigma_T \le \rho $. Since $B$ is symmetric $x^\top A^k B A^k x \le \lambda_{\max} (B) \norm{A^k x}^2$. Moreover $B$ is a rank 1 matrix and thus its sole nonzero eigenvalue is $\lambda_{\max} (B) = \norm{\beta_T - \beta_S}^2$. The aforementioned hypothesis leads to $\frac{1}{\sigma_T^2} x^\top A^k B A^k x \le  \rho^2 \norm{A^k x}^2$. Therefore we have the following lower bound $\psi_k(x)$ of $\phi_k(x)$ that can be used in practice:

\[
\psi_k(x) = \frac{\hsigma_T^2}{\hsigma_S^2} \frac{x^\top (\Sigma_T^{-1} - \alpha^2 \Omega_k \Sigma_T \Omega_k) x - \rho^2 \norm{A^k x}^2}{x^\top A^k \Sigma_S A^k x}
\]

What remains to prove is that the type I error is maximum at the frontier, i.e. where $\Delta \RR_k(x) = 0$. If $\Delta \RR_k(x) \le 0$ then:

\[
\phi_k(x) \le \frac{\hsigma_T^2 \frac{\sigma_S^2 x^\top A^k \Sigma_S^{-1} A^k x + x^\top A^k B A^k x}{\sigma_T^2} - (\hsigma_T^2 / \sigma_T^2) x^\top A^k B A^k x}{\hsigma_S^2 x^\top A^k \Sigma_S^{-1} A^k x} 
\]

with an equality on the frontier. Finally the r.h.s. can be simplified in $F= \frac{\hsigma_T^2 / \sigma_T^2}{\hsigma_S^2 / \sigma_S}$. Thus finally:

$$
\P_{\Delta \RR_k(x) \le 0}(\phi_k(x) \ge q^{1-a} ) \le \P_{F\sim \F(N_T-D,N_S-D)} \Big( F \ge q^{1-a} \Big) = \P_{\Delta \RR_k(x) = 0}(\phi_k(x) \ge q^{1-a} \big) = a
$$

which proves that the type I error is maximum at $\Delta \RR_k(x) = 0$ and that the level of the test is $a$. \\

Thus the p-value of the test relying on $\phi_k(x)$ can thus be upper bounded by $\P_{F\sim\F(N_T-D,N_S-D)} \big(F \ge \psi_k(x) \big)$, proving all the results of the theorem. 
\end{proof}

\subsection{Proof of equation (\ref{eq:KL})}

\begin{proof}
The KL divergence between two univariate gaussians directly yields:

$$
\begin{array}{lcl}
     2 D_{KL}(\NN_k || \NN_T) & = & \displaystyle\frac{x^\top V_k x}{\sigma_T^2 x^\top \Sigma_T^{-1} x} + \displaystyle\frac{(x^\top \big(\beta_T - \beta_k) \big)^2 - \sigma_T^2 x^\top \Sigma_T^{-1} x}{\sigma_T^2 x^\top \Sigma_T^{-1} x} - \ln \Big( \displaystyle\frac{x^\top V_k x}{\sigma_T^2 x^\top \Sigma_T^{-1} x} \Big) \\
     \\
     & = & \displaystyle\frac{-\Delta \RR_k(x)}{\sigma_T^2 x^\top \Sigma_T^{-1} x} - \ln \Big( \displaystyle\frac{x^\top V_k x}{\sigma_T^2 x^\top \Sigma_T^{-1} x} \Big) \quad \text{hence the result.}
\end{array}
$$
\end{proof}

%\newpage

\section{Additional experimental elements} 

The first subsection is dedicated to the calibration of the hyperparameters $k$ and $\rho$ that was discussed theoretically in section III of the paper. The second one gives more insight on the GEFCOM2012 data.

\subsection{Synthetic data - choice of the hyperparameters}

\begin{figure}[h]
    \centering
    \subfloat[][Choice of $\hat{k}$]{\includegraphics[width=0.43\linewidth]{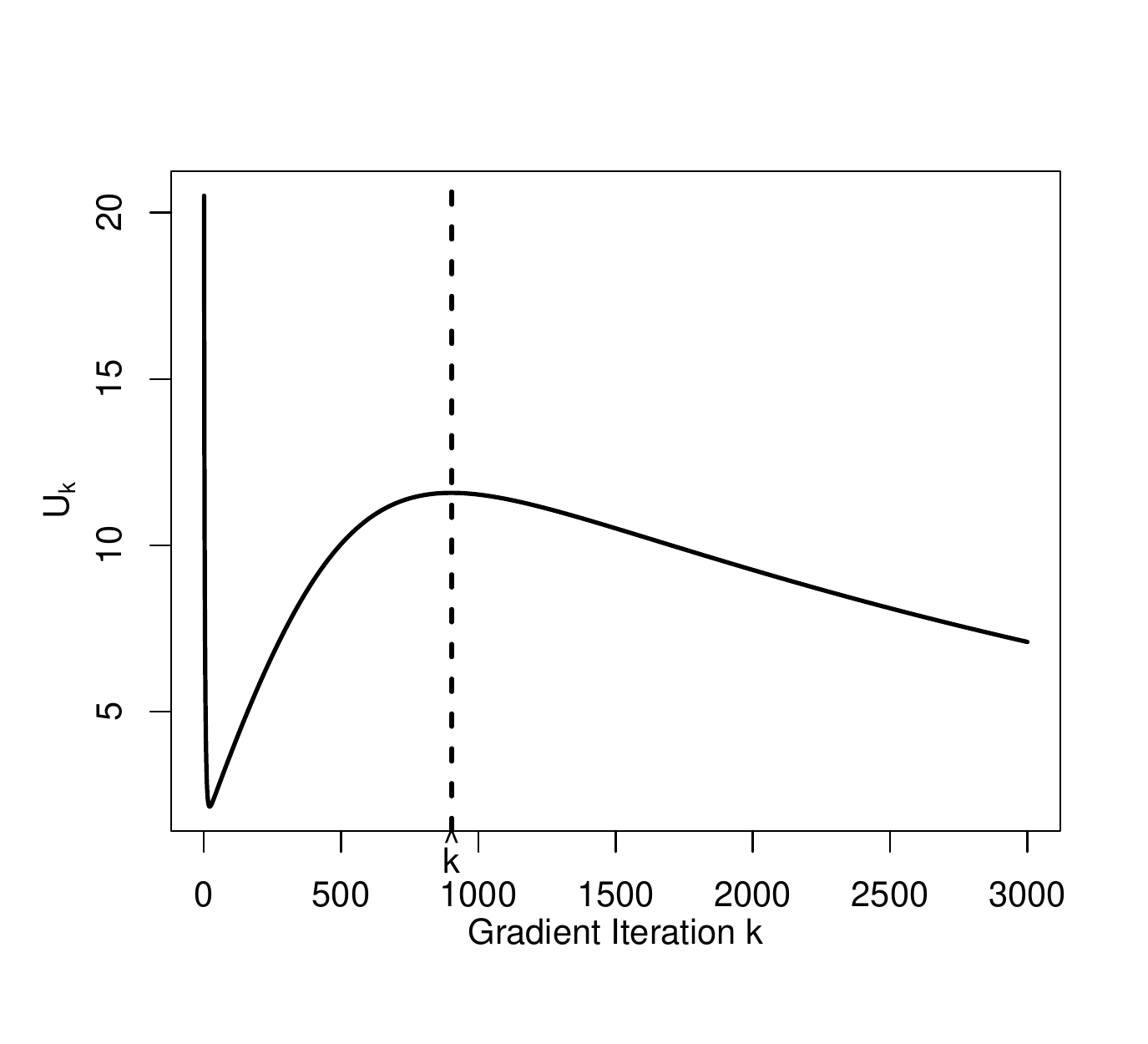}}
    \subfloat[][Choice of $\hat{\rho}$]{\includegraphics[width=0.43\linewidth]{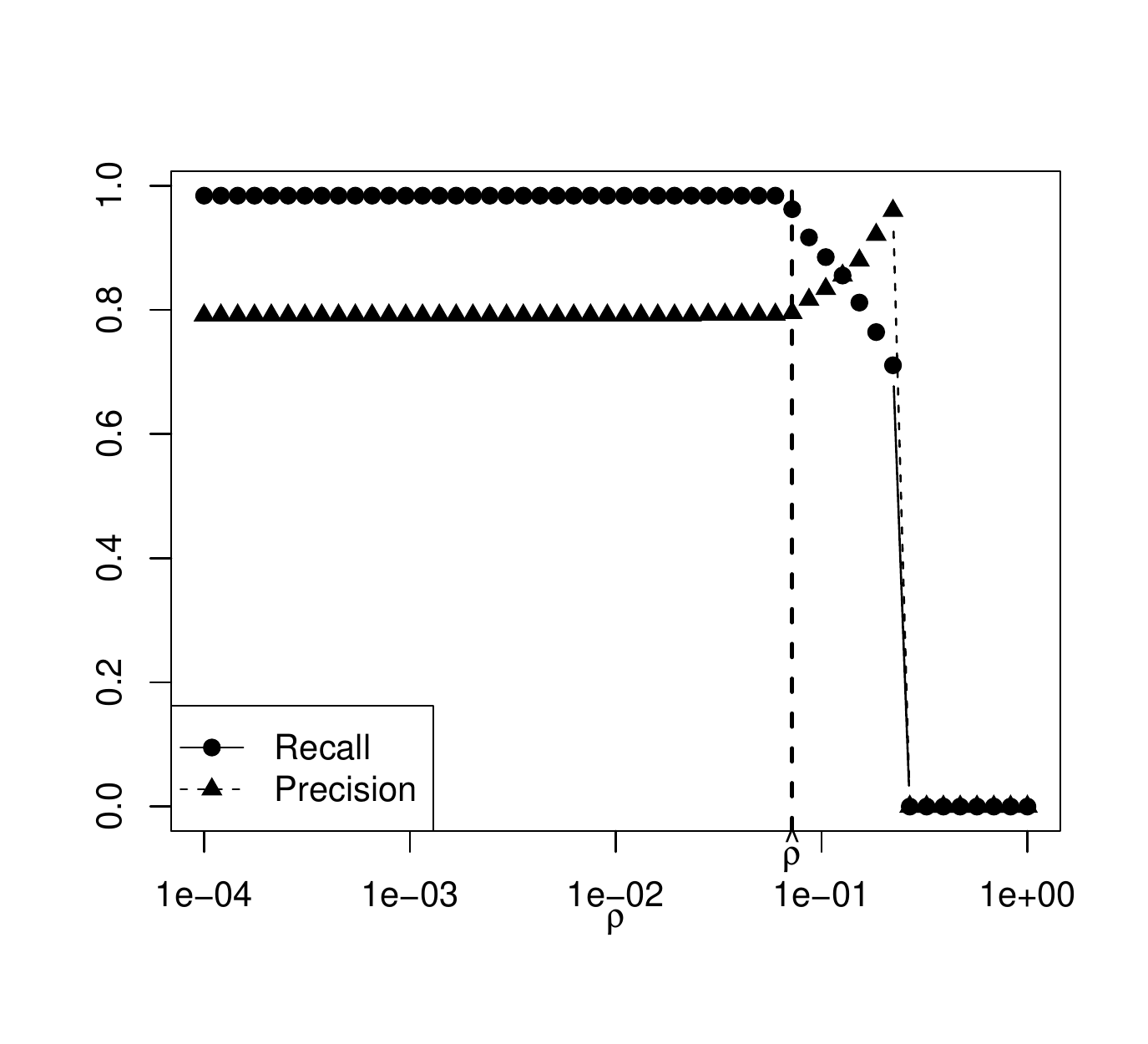}}
    \caption{Tuning of the test quantities $k$ and $\rho$}
    \label{fig:poly_choix_k_rho}
\end{figure}

The procedure to choose $k$ and $\rho$ is detailed further with figure \ref{fig:poly_choix_k_rho}. The plot of $\overline{U}_k$ used to choose $\hat{k}$ is represented in (a). Initially $\overline{U}_k$ starts from a very high value, but the number of iterations is too low. Thus the local maximum at $\hat{k}=4053$ is chosen instead. The tuning of $\rho$ is performed as following. The training samples $(x_i,y_i) \in \D_S \cup \D_T$ such that $(y_i-\hat{y}_{T,i})^2 > (y_i-\hat{y}_{k,i})^2$ are labeled $1$, and the others $0$. This can be seen as using an empirical counterpart to the gain. The test is then applied on this data for different values of $\rho$, on which we check the recall and precision for the label $1$. Finally $\hat{\rho}$ is taken to maximize the precision with a good recall or right before the latter drops, as represented in (b).
\vspace*{-1cm}

\subsection{GEFCOM2012 data representation}

\begin{figure}[h]
    \centering
    \subfloat[][Load over 2004-2005]{\includegraphics[width=0.43\linewidth]{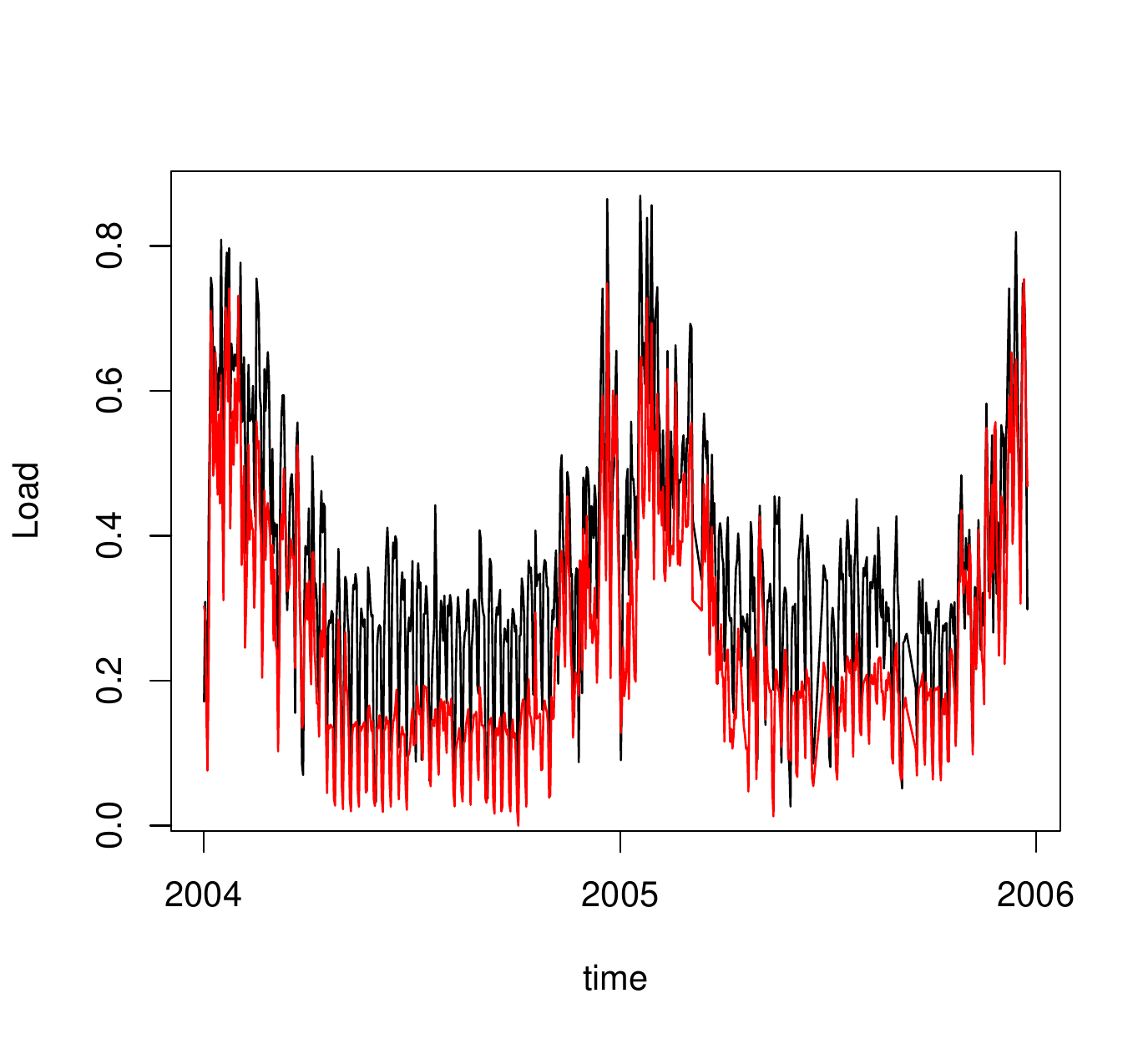}}
    \subfloat[][Daily Load]{\includegraphics[width=0.43\linewidth]{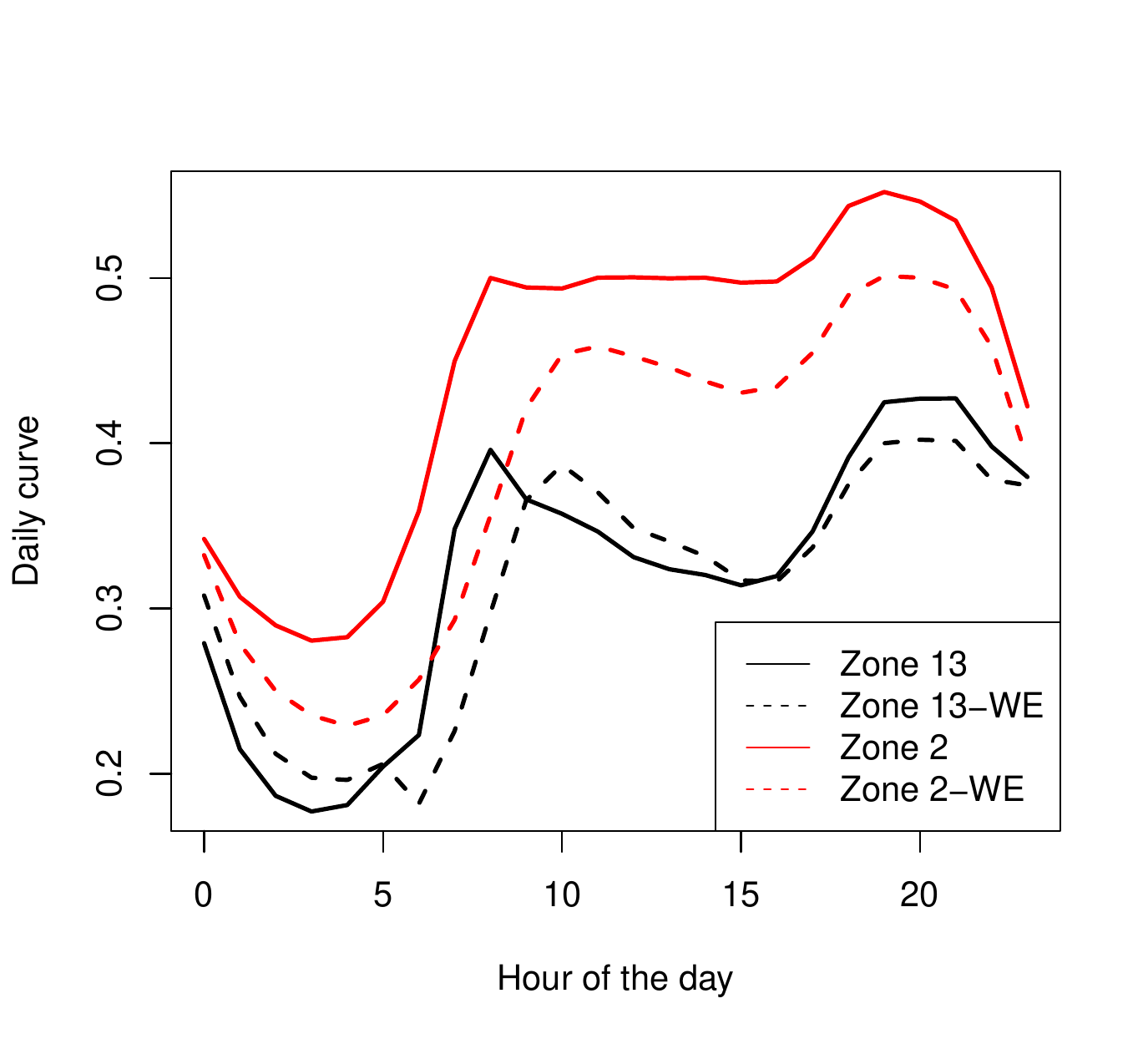}}
    \caption{Comparison of the load demand for zones 13 ($\S$) and 2 ($\T$).}
    \label{fig:comparison_loads}
\end{figure}

The electricity demand data of the GEFCOM2012 dataset is represented figure \ref{fig:comparison_loads}. Fig. (a) depicts the annual trend for the demand of the two zones at 8a.m, whereas (b) represents the average demand over a day during the week or on weekends (WE). The behavior of the two series is similar: the annual trend is the same (higher consumption in winter and lower during the summer), and the daily peaks happen around the same hour.
\vspace*{-1cm}

\bibliographystyle{apalike}% plainnat
\bibliography{biblio}

\end{document}